# Arbitrarily Accurate Analytical Approximations for the Error Function

**Roy M. Howard**

*School of Electrical Engineering, Computing and Mathematical Sciences*
*Curtin University, GPO Box U1987, Perth, 6845, Australia.*

*r.howard@curtin.edu.au*

*6 Feb 2022*

**Abstract:** A spline based integral approximation is utilized to define a sequence of approximations to the error function that converge at a significantly faster manner than the default Taylor series. The real case is considered and the approximations can be improved by utilizing the approximation $\text{erf}(x) \approx 1$ for $|x| > x_o$ and with $x_o$ optimally chosen. Two generalizations are possible; the first is based on demarcating the integration interval into $m$ equally spaced sub-intervals. The second, it based on utilizing a larger fixed sub-interval, with a known integral, and a smaller sub-interval whose integral is to be approximated. Both generalizations lead to significantly improved accuracy. Further, the initial approximations, and the approximations arising from the first generalization, can be utilized as the inputs to a custom dynamical system to establish approximations with better convergence properties.

Indicative results include those of a fourth order approximation, based on four sub-intervals, which leads to a relative error bound of $1.43 \times 10^{-7}$ over the interval $[0, \infty]$. The corresponding sixteenth order approximation achieves a relative error bound of $2.01 \times 10^{-19}$. Various approximations, that achieve the set relative error bounds of $10^{-4}$, $10^{-6}$, $10^{-10}$ and $10^{-16}$, over $[0, \infty]$, are specified.

Applications include, first, the definition of functions that are upper and lower bounds, of arbitrary accuracy, for the error function. Second, new series for the error function. Third, new sequences of approximations for $\exp(-x^2)$ which have significantly higher convergence properties that a Taylor series approximation. Fourth, the definition of a complementary demarcation function $e_C(x)$ which satisfies the constraint $e_C^2(x) + \text{erf}^2(x) = 1$. Fifth, arbitrarily accurate approximations for the power and harmonic distortion for a sinusoidal signal subject to an error function nonlinearity. Sixth, approximate expressions for the linear filtering of a step signal that is modelled by the error function.

## 1   Introduction

The error function arises in many areas of mathematics, science and scientific applications including diffusion associated with Brownian motion (Fick's second law), the heat kernel for the heat equation (e.g. Lebedev 1971), the modelling of magnetization (e.g. Fujiwara, 1980), the modelling of transitions between two levels, for example, with the modelling of smooth or soft limiters (e.g. Lee 1972) and the psychometric function (e.g. Klein 2001, Rinderknecht 2018), the modelling of amplifier non-linearities (e.g. Shi 1996 and Taggart 2005) and the modelling of rubber like materials and soft tissue (e.g. Li 2016 and Ogden 1999). It is widely used in the modelling of random phenomena as the error function defines the cumulative distribution of the Gaussian probability density function and examples include, the probability of error in signal detection, option pricing via the Black-Scholes formula etc. Many other applications exist. In general, the error function is associated with a macro description of physical phenomena and the De-Moivre Laplace theorem is illustrative of the link between fundamental outcomes and a higher level model.

The error function is defined on the complex plane according to

$$\text{erf}(z) = \frac{2}{\sqrt{\pi}} \int_\gamma e^{-\lambda^2} d\lambda, \qquad z \in \mathbf{C}, \tag{1}$$

where the path $\gamma$ is between the points zero and $z$ and is arbitrary. Associated functions are the complementary error function, the Faddeyeva function, the Voigt function and Dawson's integral (e.g. Temme 2010). The Faddeyeva function and the Voigt function, for example, have application in spectroscopy (e.g. Schrerier 1992). The error function can also be defined in terms of the spherical Bessel functions (e.g. Temme 2010, eqn. 7.6.8) and the incomplete Gamma function (e.g. Temme 2010, eqn. 7.11.1). Marsaglia (2004) provides a brief insight into the history of the error function.

For the real case, which is the case considered in this paper, the error function is defined by the integral

$$\mathrm{erf}(x) = \frac{2}{\sqrt{\pi}} \int_0^x e^{-\lambda^2} d\lambda, \qquad x \in \mathbf{R}. \tag{2}$$

The widely used, and associated, cumulative distribution function for the standard normal distribution is defined according to

$$\Phi(x) = \frac{1}{\sqrt{2\pi}} \int_{-\infty}^x \exp\left[\frac{-\lambda^2}{2}\right] d\lambda = 0.5 + 0.5\,\mathrm{erf}\left[\frac{x}{\sqrt{2}}\right]. \tag{3}$$

Being defined by an integral, which does not have an explicit analytical form, there is interest in approximations for the error function and over recent decades many approximations have been developed. Table 1 details indicative approximations for the real case and their relative errors are shown in Figure 1. Most of the approximations detailed in this table are custom and have a limited relative error bound with bounds in the range of $3.05 \times 10^{-5}$ (Sandoval-Hernandez) to $7.07 \times 10^{-3}$ (Menzel). It is preferable to have an approximation form that can be generalized to create approximations that converge to the error function. Examples include the standard Taylor series and the Bürmann series defined in Table 1.

Many of the approximations detailed in Table 1 can be improved upon by approximating the associated residual function, denoted $g$, via a Padé approximant or a basis set decomposition. Examples of some of the possible approximation forms, and the resulting residual functions, are detailed in Table 2. One example is that of a $4/2$ Padé approximant for the function $g_3$ which leads to the approximation:

$$\mathrm{erf}(x) \approx \sqrt{1 - \exp\left[-x^2 \cdot \frac{4}{\pi}\left[1 + \frac{n_1 x_1 + n_2 x_1^2 + n_3 x_1^3 + n_4 x_1^4}{1 + d_1 x_1 + d_2 x_1^2}\right]\right]}, \qquad x_1 = \frac{x}{x+1}, x \geq 0, \tag{4}$$

$$n_1 = \frac{279}{10,000,000} \qquad n_2 = \frac{-303,923}{10,000,000} \qquad n_3 = \frac{34783}{5,000,000} \qquad n_4 = \frac{40793}{10,000,000}$$
$$d_1 = \frac{-21,941,279}{10,000,000} \qquad d_2 = \frac{3,329,407}{2,500,000}. \tag{5}$$

The relative error bound in this approximation is $4.02 \times 10^{-7}$. Higher order Padé approximants can be used to generate approximations with a lower relative error bound. Matic 2018 provides a similar approximation with an absolute error of $5.79 \times 10^{-6}$.

An approximation for the error function can also be obtained by combining separate approximations, which are accurate, respectively, for $|x| \ll 1$ and $|x| \gg 1$, via a demarcation function $d$

$$\mathrm{erf}(x) = \frac{2x}{\sqrt{\pi}} \cdot d(x) + \left[1 - \frac{e^{-x^2}}{\sqrt{\pi}x}\right] \cdot [1 - d(x)], \qquad x \geq 0, \tag{6}$$

where

$$d(x) = \frac{\sqrt{\pi}x\,\mathrm{erf}(x) - \sqrt{\pi}x + e^{-x^2}}{2x^2 - \sqrt{\pi}x + e^{-x^2}}. \tag{7}$$

Naturally, an approximation for $d$ is required which requires an approximation for the error function. Unsurprisingly, the relative error in the approximation for the error function equals the relative error in the approximation utilized to approximate the error function in $d$.

Finally, efficient numerical implementation of the error function is of interest and Chevillard (2012) and De Schrijver (2018) provide results and an overview. Highly accurate piecewise approximations have long been defined, e.g. Cody 1969.

**Table 1.** Examples of published approximations for $\mathrm{erf}(x)$, $0 < x < \infty$. For the third and second last approximations the coefficient definitions are detailed in the associated reference. The stated relative error bounds arise from sampling the interval $[0, 5]$ with 10,000 uniformly spaced points.

| # | Reference | Approximation | Relative error bound |
|---|---|---|---|
| 1 | Taylor series | $T_n(x) = \dfrac{2}{\sqrt{\pi}} \cdot \left[ x - \dfrac{x^3}{3 \cdot 1} + \dfrac{x^5}{5 \cdot 2!} - \ldots + \dfrac{(-1)^{(n-1)/2} x^n}{n \cdot [(n-1)/2]!} \right]$, $n \in \{1, 3, 5, \ldots\}$ | |
| 2 | Abramowitz, 1964, p. 297, eqn. 7.1.6 | $\dfrac{2}{\sqrt{\pi}} \left[ x + \dfrac{2x^3}{1 \cdot 3} + \dfrac{2^2 x^5}{3 \cdot 5} + \dfrac{2^3 x^7}{3 \cdot 5 \cdot 7} + \ldots + \dfrac{2^n x^{2n+1}}{1 \cdot 3 \cdot 5 \cdot \ldots \cdot (2n+1)} \right] e^{-x^2}$ | |
| 3 | Abramowitz, 1964, p. 299, eqn. 7.1.26 | $1 - \left[ \dfrac{a_1}{1+px} + \dfrac{a_2}{(1+px)^2} + \dfrac{a_3}{(1+px)^3} + \dfrac{a_4}{(1+px)^4} + \dfrac{a_5}{(1+px)^5} \right] e^{-x^2}$ | $8.09 \times 10^{-6}$ |
| 4 | Menzel, 1975 Nadagopal, 2010 | $\sqrt{1 - \exp\left[\dfrac{-4x^2}{\pi}\right]}$ | $7.07 \times 10^{-3}$ |
| 5 | Bürmann series Schöpf, 2014, eqn. 33. | $\dfrac{2}{\sqrt{\pi}} \cdot \sqrt{1 - \exp(-x^2)} \cdot \left[ \dfrac{\sqrt{\pi}}{2} + \dfrac{31}{200} e^{-x^2} - \dfrac{341}{8000} e^{-2x^2} \right]$ | $3.61 \times 10^{-3}$ |
| 6 | Winitzki, 2008, eqn. 3. | $\sqrt{1 - \exp\left[-x^2 \cdot \dfrac{4/\pi + ax^2}{1 + ax^2}\right]}$, $a = \dfrac{8(\pi - 3)}{3\pi(4 - \pi)}$ | $3.50 \times 10^{-4}$ |
| 7 | Soranzo, 2012, eqn. 1. | $\sqrt{1 - \exp\left[-x^2 \cdot \dfrac{a_1 + a_2 x^2}{1 + b_2 x^2 + b_3 x^4}\right]}$ $\begin{cases} a_1 = 1.2735457 \\ a_2 = 0.1487936 \end{cases}$ $b_2 = 0.1480931 \quad b_3 = 5.160 \times 10^{-4}$ | $1.20 \times 10^{-4}$ |
| 8 | Vedder, 1987, eqn. 5. | $\tanh\left[\dfrac{167x}{148} + \dfrac{11x^3}{109}\right]$ | $4.65 \times 10^{-3}$ |
| 9 | Vazquez-Leal, 2012, eqn. 3.1. | $\tanh\left[\dfrac{39x}{2\sqrt{\pi}} - \dfrac{111}{2} \cdot \mathrm{atan}\left[\dfrac{35x}{111\sqrt{\pi}}\right]\right]$ | $1.88 \times 10^{-4}$ |
| 10 | Sandoval-Hernandez, 2019, eqn. 23. | $\dfrac{2}{1 + \exp[\alpha_1 x + \alpha_3 x^3 + \alpha_5 x^5 + \alpha_7 x^7 + \alpha_9 x^9]} - 1$ | $3.05 \times 10^{-5}$ |
| 11 | Abrarov, 2013, eqn. 16. | $1 - e^{-x^2}\left[\dfrac{1 - e^{-\tau_m x}}{\tau_m x} + \dfrac{\tau_m^2 x}{\sqrt{\pi}} \cdot \sum_{n=1}^{N} \dfrac{a_n[1 - (-1)^n e^{-\tau_m x}]}{n^2 \pi^2 + \tau_m^2 x^2}\right]$ $a_n = \dfrac{2\sqrt{\pi}}{\tau_m} \cdot \exp[-n^2 \pi^2 / \tau_m^2]$, $\tau_m = 12$ | $3.27 \times 10^{-3}$ ($N = 6$) |

**Table 2.** Residual functions associated with approximations for $\mathrm{erf}(x)$, $0 < x < \infty$.

| # | Error function | Residual Function |
|---|---|---|
| 1 | $\tanh\left[\dfrac{2x}{\sqrt{\pi}}\right] + g_1(x)$ | $g_1(x) = \mathrm{erf}(x) - \tanh\left[\dfrac{2x}{\sqrt{\pi}}\right]$ |
| 2 | $\tanh\left[\dfrac{2x}{\sqrt{\pi}}[1 + g_2(x)]\right]$ | $g_2(x) = \dfrac{\sqrt{\pi}}{2x} \cdot \mathrm{atanh}[\mathrm{erf}(x)] - 1$ |
| 3 | $\sqrt{1 - \exp\left[-x^2 \cdot \dfrac{4}{\pi}[1 + g_3(x)]\right]}$ | $g_3(x) = \dfrac{-\pi}{4x^2} \cdot \ln[1 - \mathrm{erf}(x)^2] - 1$ |

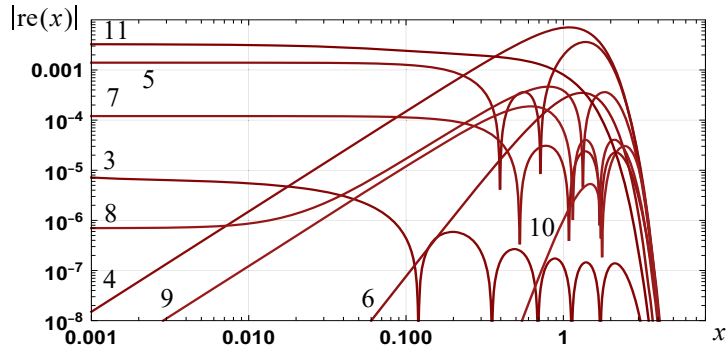

**Figure 1.** Graph of the magnitude of the relative error in the approximations, detailed in Table 1, for erf($x$).

The two point spline based approximations for functions and integrals, detailed in Howard 2019, have recently been applied to find arbitrarily accurate approximations for the hyperbolic tangent function (Howard 2021). In this paper, the general two point spline approximation form is applied to define a sequence of convergent approximations for the error function. The basic form of the approximation of order $n$, denoted $f_n$, is

$$\text{erf}(x) \approx f_n(x) = p_{n,0}(x) + p_{n,1}(x)e^{-x^2} \tag{8}$$

where $p_{n,1}$ is a polynomial of order $n$ and $p_{n,0}$ is a polynomial of order less than $n$. Convergence of the sequence of approximations $f_1, f_2, \ldots$ to erf($x$) is shown and the convergence is significantly better than the default Taylor series. For example, the second order approximation

$$f_2(x) = \frac{x}{\sqrt{\pi}} \cdot \left[1 - \frac{x^2}{30}\right] + \frac{x}{\sqrt{\pi}} \cdot \left[1 + \frac{11x^2}{30} + \frac{x^4}{15}\right]e^{-x^2} \tag{9}$$

yields as relative error bound of 0.056 over the interval $[0, 2]$ which is better that a fifteenth order Taylor series approximation. The approximations can be improved by utilizing the approximation erf($x$) ≈ 1 for $|x| \gg 1$ and thereby establishing approximations with a set relative error bound over the interval $[0, \infty)$.

Two generalizations are detailed. The first is of the form

$$\text{erf}(x) \approx p_0(x) + p_1(x)e^{-k_1 x^2} + \ldots + p_m(x)e^{-k_m x^2} \tag{10}$$

and is based on utilizing approximations associated with $m$ equally spaced sub-intervals of the interval $[0, x]$. The second, is based on utilizing a fixed sub-interval within $[0, x]$ and then approximating the error function over the remainder of the interval. Both generalizations lead to significantly improved accuracy. For example, a fourth order approximation based on four variable sub-intervals, when used with the approximation erf($x$) ≈ 1 for $x \gg 1$, has a relative error bound of $1.43 \times 10^{-7}$ over the interval $[0, \infty]$. The corresponding sixteenth order approximation has a relative error bound of $2.01 \times 10^{-19}$. Finally, by utilizing the solutions of a custom dynamical system, approximations with better convergence properties can be established.

Applications of the proposed approximations for the error function include, first, the definition of functions that are upper and lower bounds, of arbitrary accuracy, for the error function. Second, new series for the error function. Third, new sequences of approximations for $\exp(-x^2)$ which have significantly higher convergence properties that a Taylor series approximation. Fourth, the definition of a complementary demarcation function $e_C(x)$ which satisfies the constraint $e_C^2(x) + \text{erf}^2(x) = 1$. Fifth, arbitrarily accurate approximations for the power and harmonic distortion for a sinusoidal signal subject to a error function nonlinearity. Sixth, approximate expressions for the linear filtering of a step signal modelled by the error function.

Section 2 details the spline based approximation for the error function and its convergence. Improved approximations, obtained by utilizing the nature of the error function for large arguments, are detailed in section 3. Two generalizations, with potential for lower relative error bounds, are detailed in sections 4 and 5. Section 6 details how the initial approximations, and the approximations arising from the first generalization, can be utilized as the inputs to a custom dynamical system to establish approximations with better convergence properties. Applications are specified in section 7. Conclusions are stated in section 8.

### 1.1   Notes and Notation

As $\text{erf}(-x) = -\text{erf}(x)$ it is sufficient to consider approximations for the interval $[0, \infty)$.

For a function $f$ defined over the interval $[\alpha, \beta]$, an approximating function $f_A$ has a relative error, at a point $x_1$, defined according to $\text{re}(x_1) = 1 - f_A(x_1)/f(x_1)$. The relative error bound for the approximating function over the interval $[\alpha, \beta]$ is defined according to

$$\text{re}_B = \max\{|\text{re}(x_1)|: x_1 \in [\alpha, \beta]\}. \tag{11}$$

The notation $f^{(k)}(x) = \dfrac{d^k}{dx^k} f(x)$ is used. The symbol $u$ denotes the unit step function.

Mathematica has been used to facilitate analysis and to obtain numerical results.

### 1.2   Background Results

The following result underpins the bounds proposed for the error function:

#### Lemma 1    Upper and Lower Functional Bounds

A positive approximation $f_A$ to a positive function $f$ over the interval $[\alpha, \beta]$, with a relative error bound

$$-\varepsilon_B < 1 - \frac{f_A(x)}{f(x)} < \varepsilon_B, \qquad x \in [\alpha, \beta], \varepsilon_B > 0, \tag{12}$$

leads to the following upper and lower bounded functions:

$$\frac{f_A(x)}{1 + \varepsilon_B} < f(x) < \frac{f_A(x)}{1 - \varepsilon_B}, \qquad x \in [\alpha, \beta]. \tag{13}$$

The relative error bounds, over the interval $[\alpha, \beta]$, for the upper and lower bounded functions, respectively, are:

$$\frac{2\varepsilon_B}{1 + \varepsilon_B} \qquad \frac{2\varepsilon_B}{1 - \varepsilon_B} \tag{14}$$

#### Proof

The definition of the relative error bound, as specified by Equation 12, leads to

$$1 - \varepsilon_B < \frac{f_A(x)}{f(x)} < 1 + \varepsilon_B \tag{15}$$

which implies

$$\frac{1 - \varepsilon_B}{1 + \varepsilon_B} < \frac{f_A(x)/(1 + \varepsilon_B)}{f(x)} < 1, \qquad 1 < \frac{f_A(x)/(1 - \varepsilon_B)}{f(x)} < \frac{1 + \varepsilon_B}{1 - \varepsilon_B}, \tag{16}$$

and the relative error bounds:

$$0 < 1 - \frac{f_A(x)/(1+\varepsilon_B)}{f(x)} < 1 - \frac{1-\varepsilon_B}{1+\varepsilon_B} = \frac{2\varepsilon_B}{1+\varepsilon_B}, \qquad 1 - \frac{1+\varepsilon_B}{1-\varepsilon_B} = \frac{-2\varepsilon_B}{1-\varepsilon_B} < 1 - \frac{f_A(x)/(1-\varepsilon_B)}{f(x)} < 0. \tag{17}$$

#### 1.2.1   Convergent Integral Approximations

One application of the proposed approximations for the error function requires knowledge of when function convergence implies convergence of associated integrals.

#### Lemma 2    Convergent Integral Approximation

If a sequence of functions $f_1, f_2, \ldots$ converges, over the interval $[0, x]$, to a bounded and integrable function $f$, then point-wise convergence is sufficient for the associated integrals to be convergent, i.e. for

$$\lim_{n \to \infty} \int_0^x f_n(\lambda) d\lambda = \int_0^x f(\lambda) d\lambda. \tag{18}$$

**Proof**

The required result follows if it is possible to interchange the order of limit and integration, i.e.

$$\lim_{n \to \infty} \int_0^x f_n(\lambda) d\lambda = \int_0^x \lim_{n \to \infty} f_n(\lambda) d\lambda = \int_0^x f(\lambda) d\lambda. \tag{19}$$

Standard conditions for when the interchange is valid are specified by the monotone and dominated convergence theorems (e.g. Champeney 1987, p. 26). Sufficient conditions for a valid interchange include point-wise function convergence, and for $f$ to be integrable and bounded.

## 2   Spline Based Approximations for Error Function

### 2.1   Spline Approximation for Error Function

The following $n$th order, two point spline based, approximation for an integral has been detailed in Howard 2019 (eqn. 48) for a function $f$ that is at least $(n+1)$th order differentiable:

$$\int_\alpha^x f(\lambda) d\lambda = \sum_{k=0}^n c_{n,k}(x-\alpha)^{k+1}[f^{(k)}(\alpha) + (-1)^k f^{(k)}(x)] + R_n(\alpha, x) \tag{20}$$

$$n \in \{0, 1, 2, \ldots\}$$

where

$$c_{n,k} = \frac{n!}{(n-k)!(k+1)!} \cdot \frac{(2n+1-k)!}{2(2n+1)!}, \quad k \in \{0, 1, \ldots, n\}, \tag{21}$$

$$R_n^{(1)}(\alpha, x) = -\sum_{k=0}^n c_{n,k}(k+1)(x-\alpha)^k \left[ f^{(k)}(\alpha) + (-1)^{k+1} \cdot \frac{n+1}{n-k+1} \cdot f^{(k)}(x) \right] + \tag{22}$$

$$c_{n,n}(-1)^{n+1}(x-\alpha)^{n+1} f^{(n+1)}(x).$$

Direct application of this result to the integral defining the error function leads to the following result:

**Theorem 2.1   Spline Based Integral Approximation for Error Function**

The error function can be defined according to

$$\text{erf}(x) = f_n(x) + \varepsilon_n(x) \tag{23}$$

where $f_n$ is the $n$th order spline based integral approximation defined according to

$$f_n(x) = \frac{2}{\sqrt{\pi}} \cdot \sum_{k=0}^n c_{n,k} x^{k+1} \left[ p(k, 0) + (-1)^k p(k, x) e^{-x^2} \right] \tag{24}$$

and $\varepsilon_n(x)$ is the associated residual function whose derivative is

$$\varepsilon_n^{(1)}(x) = \frac{2e^{-x^2}}{\sqrt{\pi}} - \frac{2}{\sqrt{\pi}} \cdot \sum_{k=0}^n c_{n,k}(k+1)x^k p(k, 0) -$$

$$\frac{2e^{-x^2}}{\sqrt{\pi}} \cdot \sum_{k=0}^n c_{n,k} x^k (-1)^k [(k+1-2x^2) p(k, x) + x p^{(1)}(k, x)] \tag{25}$$

In these equations

$$p(k, x) = p^{(1)}(k-1, x) - 2xp[k-1, x], \qquad p(0, x) = 1. \tag{26}$$

A more general approximation is

$$\operatorname{erf}(x) - \operatorname{erf}(\alpha) = \frac{2}{\sqrt{\pi}} \cdot \sum_{k=0}^{n} c_{n,k}(x-\alpha)^{k+1} \left[ p(k, \alpha)e^{-\alpha^2} + (-1)^k p(k, x)e^{-x^2} \right] + \varepsilon_n(\alpha, x)$$

$$= f_n(\alpha, x) + \varepsilon_n(\alpha, x). \tag{27}$$

**Proof**

The proof is detailed in Appendix 1.

### 2.1.1   Note

The polynomial function $p(k, x)$ is equivalently defined by the $k$th order Hermite polynomial function (Abramowitz, 1964, p. 775, equation 23.3.10) and an explicit form is

$$p(k, x) = \sum_{i=0}^{\lfloor k/2 \rfloor} \frac{(-1)^{i+k} k!}{i!(k-2i)!} \cdot 2^{k-2i} x^{k-2i}. \tag{28}$$

A change of variable $r = k - 2i$, and noting that $i \in \{0, 1, ..., \lfloor k/2 \rfloor\}$ implies $r \in \{k, k-2, ..., k - 2\lfloor k/2 \rfloor\}$, leads to the alternative form

$$p(k, x) = \sum_{r = k - 2\lfloor k/2 \rfloor}^{k} d_{k,r} x^r, \qquad d_{k,r} = \frac{(-1)^{(3k-r)/2} \left[ \frac{1 + (-1)^{r - [k - 2\lfloor k/2 \rfloor]}}{2} \right] k! 2^r}{[(k-r)/2]! r!}. \tag{29}$$

Substitution of this form into Equation 24 leads to the direct polynomial form for the $n$th order approximation to the error function:

$$f_n(x) = \frac{2}{\sqrt{\pi}} \cdot \sum_{k=0}^{n} c_{n,k} p(k, 0) x^{k+1} + \frac{2e^{-x^2}}{\sqrt{\pi}} \cdot \sum_{k=0}^{n} \left[ (-1)^k c_{n,k} \sum_{r = k - 2\lfloor k/2 \rfloor}^{k} d_{k,r} x^{r+k+1} \right]. \tag{30}$$

### 2.1.2   Approximations

The polynomial functions $p$, as defined by Equation 26 or Equation 28, have the explicit forms:

$$p(0, x) = 1, \qquad p(1, x) = -2x, \qquad p(2, x) = -2[1 - 2x^2],$$
$$p(3, x) = 12x[1 - 2x^2/3], \qquad p(4, x) = 12[1 - 4x^2 + 4x^4/3], \qquad \ldots \tag{31}$$

Approximations for the error function, as defined by Equation 24 and for orders zero to five, are:

$$f_0(x) = \frac{x}{\sqrt{\pi}} + \frac{x}{\sqrt{\pi}} \cdot e^{-x^2} \tag{32}$$

$$f_1(x) = \frac{x}{\sqrt{\pi}} + \frac{x}{\sqrt{\pi}} \left[ 1 + \frac{x^2}{3} \right] e^{-x^2} \tag{33}$$

$$f_2(x) = \frac{x}{\sqrt{\pi}} \left[ 1 - \frac{x^2}{30} \right] + \frac{x}{\sqrt{\pi}} \left[ 1 + \frac{11x^2}{30} + \frac{x^4}{15} \right] e^{-x^2} \tag{34}$$

$$f_3(x) = \frac{x}{\sqrt{\pi}} \left[ 1 - \frac{x^2}{21} \right] + \frac{x}{\sqrt{\pi}} \left[ 1 + \frac{8x^2}{21} + \frac{17x^4}{210} + \frac{x^6}{105} \right] e^{-x^2} \tag{35}$$

$$f_4(x) = \frac{x}{\sqrt{\pi}}\left[1 - \frac{x^2}{18} + \frac{x^4}{1260}\right] + \frac{x}{\sqrt{\pi}}\left[1 + \frac{7x^2}{18} + \frac{37x^4}{420} + \frac{4x^6}{315} + \frac{x^8}{945}\right]e^{-x^2} \tag{36}$$

$$f_5(x) = \frac{x}{\sqrt{\pi}}\left[1 - \frac{2x^2}{33} + \frac{x^4}{660}\right] + \frac{x}{\sqrt{\pi}}\left[1 + \frac{13x^2}{33} + \frac{61x^4}{660} + \frac{67x^6}{4620} + \frac{16x^8}{10395} + \frac{x^{10}}{10395}\right]e^{-x^2} \tag{37}$$

### 2.2    Results

The relative error in the zero to tenth order spline based series approximations, along with the relative error in Taylor series approximations of order one to fifteen, are detailed in Figure 2. The clear superiority, is terms of convergence, of the spline based series, relative to the Taylor series, is evident. The relative error in the spline approximations, of orders 16, 20, 24, 28 and 32, are shown in Figure 3.

The Mathematica code underpinning the results shown in Figure 2, is detailed in Appendix 2. Such code is indicative of the code underpinning the results detailed in the paper.

### 2.3    Approximation for Large Arguments

Zero and first order approximations for the error function, and for the case of $x \gg 1$, are

$$\text{erf}(x) \approx 1, \qquad \text{erf}(x) \approx 1 - \frac{e^{-x^2}}{\sqrt{\pi}x}. \tag{38}$$

The relative errors in such approximations, respectively, are

$$\text{re}(x) = 1 - \frac{1}{\text{erf}(x)}, \qquad \text{re}(x) \approx 1 - \frac{1 - e^{-x^2}/\sqrt{\pi}x}{\text{erf}(x)}, \tag{39}$$

and their graphs are shown in Figure 3.

### 2.4    Convergence

To prove convergence of the sequence of functions $f_0, f_1, f_2, \ldots$, defined by Theorem 2.1, to the error function, it is sufficient to prove that the corresponding sequence of residual functions $\varepsilon_0, \varepsilon_1, \varepsilon_2, \ldots$ converge to zero. This can be shown by considering the derivatives of the residual functions defined by Equation 25. The derivatives of the residual functions of orders zero, one and two, respectively, are:

$$\varepsilon_0^{(1)}(x) = \frac{1}{\sqrt{\pi}}[1 + 2x^2]e^{-x^2} - \frac{1}{\sqrt{\pi}} \tag{40}$$

$$\varepsilon_1^{(1)}(x) = \frac{1}{\sqrt{\pi}}\left[1 + x^2 + \frac{2x^4}{3}\right]e^{-x^2} - \frac{1}{\sqrt{\pi}} \tag{41}$$

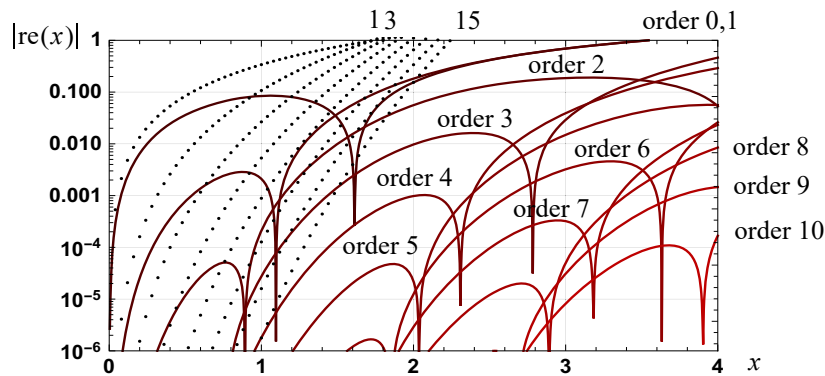

**Figure 2.**  Graph of the magnitude of the relative errors in approximations to $\text{erf}(x)$: zero to tenth order integral spline based series and first, third, ..., fifteenth order Taylor series (dotted).

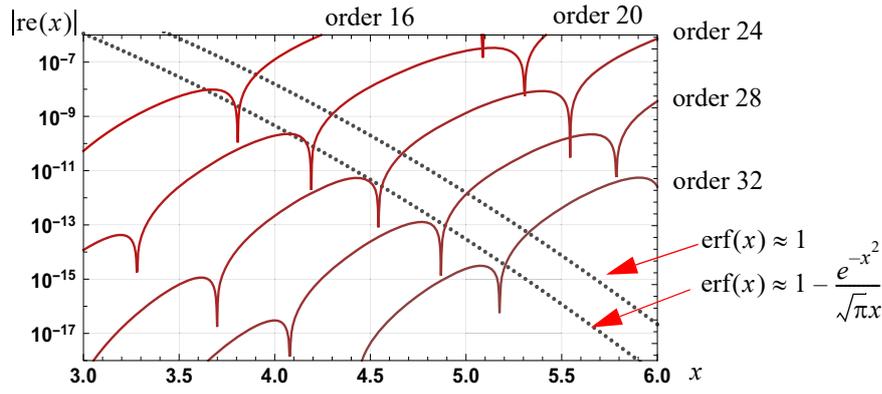

**Figure 3.** Graph of the magnitude of the relative errors associated with the approximation $\text{erf}(x) \approx 1$ and $\text{erf}(x) \approx 1 - \exp(-x^2)/\sqrt{\pi}x$ along with the relative error in spline approximations of orders 16, 20, 24, 28 and 32.

$$\varepsilon_2^{(1)}(x) = \frac{1}{\sqrt{\pi}}\left[1 + \frac{9x^2}{10} + \frac{2x^4}{5} + \frac{2x^6}{15}\right]e^{-x^2} - \frac{1}{\sqrt{\pi}}\left[1 - \frac{x^2}{10}\right]. \tag{42}$$

### Theorem 2.2  Convergence of Spline Based Approximations

For all fixed values of $x$, the derivatives of the residual functions converge to zero as the order increases, i.e. for all fixed values of $x$ it is the case that

$$\lim_{n \to \infty} \varepsilon_n^{(1)}(x) = 0, \qquad x > 0. \tag{43}$$

This is sufficient for the convergence of the residual functions, i.e. $\lim_{n \to \infty} \varepsilon_n(x) = 0$, $x > 0$, and, hence, for $x$ fixed:

$$\lim_{n \to \infty} f_n(x) = \text{erf}(x), \qquad x > 0. \tag{44}$$

The convergence is non-uniform.

### Proof

The proof is detailed in Appendix 3.

### 2.5  Improved Approximation via Iteration

Consider the general result

$$\int_0^x \text{erf}(\lambda)d\lambda = \frac{1}{\sqrt{\pi}}\left[-1 + e^{-x^2}\right] + x\text{erf}(x). \tag{45}$$

By using the approximations, $f_n$, as defined in Theorem 2.1, in the integral, improved approximations for the error function can be defined.

### Theorem 2.3  Improved Approximation via Iteration

An improved approximation, of order $n$, for the error function is

$$F_n(x) = \frac{1}{\sqrt{\pi}x}\left[1 - e^{-x^2}\right] + \frac{2}{\sqrt{\pi}} \cdot \sum_{k=0}^{n} \frac{c_{n,k}p(k,0)}{k+2} \cdot x^{k+1} + \tag{46}$$

$$\frac{2}{x\sqrt{\pi}} \cdot \sum_{k=0}^{n}\left[(-1)^k c_{n,k} \sum_{r=k-2\lfloor k/2\rfloor}^{k} \frac{d_{k,r}}{2} \cdot \left[\frac{r+k}{2}\right]! \cdot \left[1 - e^{-x^2} \cdot \sum_{i=0}^{\frac{r+k}{2}} \frac{x^{2i}}{i!}\right]\right]$$

where $c_{n,k}$ and $d_{k,r}$ are defined, respectively, by Equation 21 and Equation 29.

**Proof**

From Equation 45 it follows that

$$\text{erf}(x) \approx \frac{1}{\sqrt{\pi}x}\left[1 - e^{-x^2}\right] + \frac{1}{x} \cdot \int_0^x f_n(\lambda)d\lambda. \tag{47}$$

As

$$\int_0^x x^u e^{-\lambda^2} d\lambda = \frac{1}{2}\left[\frac{u-1}{2}\right]! \cdot \left[1 - e^{-x^2} \cdot \sum_{i=0}^{(u-1)/2} \frac{x^{2i}}{i!}\right], \qquad u \in \{1, 3, 5, \ldots\}, \tag{48}$$

it follows, from the form for $f_n$ detailed in Equation 30, that

$$\text{erf}(x) \approx \frac{1}{\sqrt{\pi}x}\left[1 - e^{-x^2}\right] + \frac{2}{\sqrt{\pi}} \cdot \sum_{k=0}^{n} \frac{c_{n,k} p(k, 0)}{k+2} \cdot x^{k+1} +$$

$$\frac{2}{x\sqrt{\pi}} \cdot \sum_{k=0}^{n}\left[(-1)^k c_{n,k} \sum_{r = k-2\lfloor k/2 \rfloor}^{k} \frac{d_{k,r}}{2} \cdot \left[\frac{r+k}{2}\right]! \cdot \left[1 - e^{-x^2} \cdot \sum_{i=0}^{\frac{r+k}{2}} \frac{x^{2i}}{i!}\right]\right] \tag{49}$$

which is the required result.

### 2.5.1   Explicit Approximations

Approximations to the error function, of orders zero to five, are:

$$f_0(x) = \frac{3}{2\sqrt{\pi}x}\left[1 + \frac{x^2}{3}\right] - \frac{3e^{-x^2}}{2\sqrt{\pi}x} \tag{50}$$

$$f_1(x) = \frac{5}{3\sqrt{\pi}x}\left[1 + \frac{3x^2}{10}\right] - \frac{5}{3\sqrt{\pi}x}\left[1 + \frac{x^2}{10}\right]e^{-x^2} \tag{51}$$

$$f_2(x) = \frac{7}{4\sqrt{\pi}x}\left[1 + \frac{2x^2}{7} - \frac{x^4}{210}\right] - \frac{7}{4\sqrt{\pi}x}\left[1 + \frac{x^2}{7} + \frac{2x^4}{105}\right]e^{-x^2} \tag{52}$$

$$f_3(x) = \frac{9}{5\sqrt{\pi}x}\left[1 + \frac{5x^2}{18} - \frac{5x^4}{756}\right] - \frac{9}{5\sqrt{\pi}x}\left[1 + \frac{x^2}{6} + \frac{23x^4}{756} + \frac{x^6}{378}\right]e^{-x^2} \tag{53}$$

$$f_4(x) = \frac{11}{6\sqrt{\pi}x}\left[1 + \frac{3x^2}{11} - \frac{x^4}{132} + \frac{x^6}{13860}\right] - \frac{11}{6\sqrt{\pi}x}\left[1 + \frac{2x^2}{11} + \frac{5x^4}{132} + \frac{16x^6}{3465} + \frac{x^8}{3465}\right]e^{-x^2} \tag{54}$$

$$f_5(x) = \frac{13}{7\sqrt{\pi}x}\left[1 + \frac{7x^2}{26} - \frac{7x^4}{858} + \frac{7x^6}{51480}\right] - \frac{13}{7\sqrt{\pi}x}\left[1 + \frac{5x^2}{26} + \frac{37x^4}{858} + \frac{313x^6}{51480} + \frac{7x^8}{12870} + \frac{x^{10}}{38610}\right]e^{-x^2} \tag{55}$$

Note that integration of these expressions leads to functions defined, in part, on the Gamma function which is an integral. This makes further iteration impractical.

### 2.5.2   Results

The relative errors in even order approximations, of orders zero to ten, are shown in Figure 4. A comparison of the results detailed in Figure 2 and Figure 4 show the clear improvement in the approximations specified by Equation 46.

## 3  Improved Approximations

### 3.1  Improved Approximation for Error Function

An improved approximation for the error function can be achieved by noting, as illustrated in Figure 3, that the approximation $\text{erf}(x) \approx 1$ is increasingly accurate for the case of $x \gg 1$ and for $x$ increasing. By switching at a suitable point $x_o$, as illustrated in Figure 5, from a spline based approximation to the approximation $\text{erf}(x) \approx 1$, an improved approximation is achieved. Naturally, it is possible to switch to the approximation $\text{erf}(x) \approx 1 - e^{-x^2}/\sqrt{\pi}x$, or higher order approximations, in a similar manner.

**Theorem 3.1  Improved Approximation for Error Function**

An improved approximation for the error function, based on a $n$th order spline approximation detailed in Theorem 2.1 or Theorem 2.3, and consistent with the illustration shown in Figure 5, is

$$\text{erf}(x) \approx f_n(x)u[x_o(n) - x] + [1 - u[x_o(n) - x]]$$
$$\text{erf}(x) \approx F_n(x)u[x_o(n) - x] + [1 - u[x_o(n) - x]] \tag{56}$$

where the transition points, respectively, are defined according to

$$x_o(n) = x: \left|1 - \frac{1}{\text{erf}(x)}\right| = \left|1 - \frac{f_n(x)}{\text{erf}(x)}\right|$$
$$x_o(n) = x: \left|1 - \frac{1}{\text{erf}(x)}\right| = \left|1 - \frac{F_n(x)}{\text{erf}(x)}\right| \tag{57}$$

**Proof**

The improved approximation results follow from optimally switching, as illustrated in Figure 5 and at the point specified by Equation 57, to the approximation $\text{erf}(x) \approx 1$ which has a lower relative error magnitude.

#### 3.1.1  Transition Points and Relative Error Bounds

The transition points, for various orders of spline approximation, are specified in Table 3. The relationship between the transition point and order is shown in Figure 6 for the case of the approximations $f_n$. This relationship can be approximated, with a second order polynomial, according to

$$x_o(n) = 1.3607 + 0.20511n - 0.002932n^2, \qquad 0 \le n \le 24. \tag{58}$$

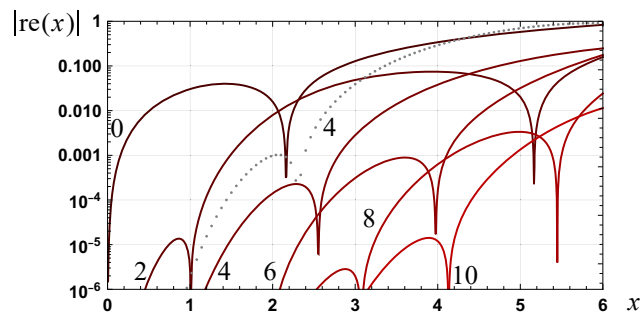

**Figure 4.** Graph of the magnitude of the relative errors in the approximations to $\text{erf}(x)$, of even orders, as specified by Equation 46. The dotted results are for the fourth order approximation specified by Theorem 2.1 (Equation 36).

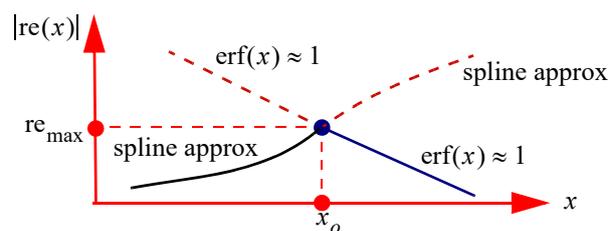

**Figure 5.** Illustration of the crossover point where the magnitude of the relative error in the approximation $\text{erf}(x) \approx 1$ equals the magnitude of the relative error in a set order spline approximation.

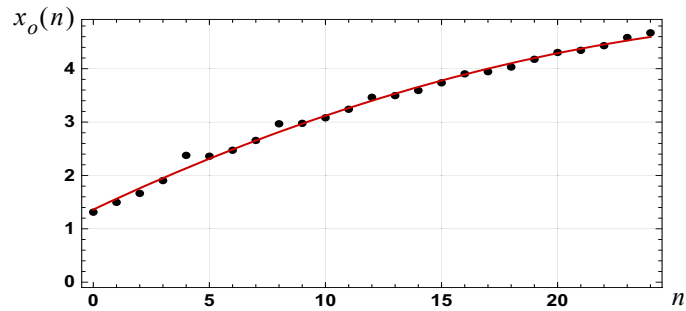

**Figure 6.** Graph of the relationship between the optimum transition point $x_o(n)$, as defined by Equation 57, and the order of the spline approximation.

However, as small variations in $x_o(n)$ can lead to significant changes in the maximum relative error in the approximation for the error function, precise values for $x_o(n)$ are preferable.

The graphs of the relative errors in the approximations $f_n$ to $\text{erf}(x)$, as specified by Equation 56, are shown in Figure 7 for orders 2, 4, 6, …, 20. The relative error bounds that can be achieved, over the interval $[0, \infty]$, using the optimally chosen transition points, are detailed in Table 3.

**Table 3.** The transition points, $x_o$, and the resulting relative error bounds for the spline based approximations specified by Equation 56. The transition points are based on sampling the interval $[0, 5]$ with 10000 points.

| Approx. order: $n$ | Transition point for $f_n$ | Relative error bound for $f_n$ | Transition point for $F_n$ | Relative error bound for $F_n$ |
|---|---|---|---|---|
| 0 | 1.3085 | 0.0851 | 1.465 | 0.0400 |
| 1 | 1.492 | 0.0362 | 1.769 | 0.0126 |
| 2 | 1.658 | $1.95 \times 10^{-2}$ | 1.929 | $6.42 \times 10^{-3}$ |
| 3 | 1.8975 | $7.36 \times 10^{-3}$ | 2.1725 | $2.13 \times 10^{-3}$ |
| 4 | 2.3715 | $1.03 \times 10^{-3}$ | 2.6305 | $2.28 \times 10^{-4}$ |
| 6 | 2.4715 | $4.75 \times 10^{-4}$ | 2.73 | $1.13 \times 10^{-4}$ |
| 8 | 2.963 | $2.79 \times 10^{-5}$ | 3.1855 | $6.69 \times 10^{-6}$ |
| 10 | 3.0785 | $1.35 \times 10^{-5}$ | 3.324 | $2.59 \times 10^{-6}$ |
| 12 | 3.4625 | $9.78 \times 10^{-7}$ | 3.67 | $2.12 \times 10^{-7}$ |
| 14 | 3.5845 | $4.00 \times 10^{-7}$ | 3.8205 | $6.57 \times 10^{-8}$ |
| 16 | 3.9025 | $3.44 \times 10^{-8}$ | 4.101 | $6.66 \times 10^{-9}$ |
| 18 | 4.0285 | $1.22 \times 10^{-8}$ | 4.257 | $1.75 \times 10^{-9}$ |
| 20 | 4.300 | $1.20 \times 10^{-9}$ | 4.493 | $2.11 \times 10^{-10}$ |
| 22 | 4.429 | $3.76 \times 10^{-10}$ | 4.652 | $4.75 \times 10^{-11}$ |
| 24 | 4.6655 | $4.18 \times 10^{-11}$ | 4.854 | $6.70 \times 10^{-12}$ |

### 3.2   Improved Approximation for Taylor Series

The approximation of $\text{erf}(x) \approx 1$ for $x \gg 1$ can be utilized to improve the relative error bound for a Taylor series approximation according to

$$\text{erf}(x) \approx T_n(x)u[x_o(n) - x] + [1 - u[x_o(n) - x]] \tag{59}$$

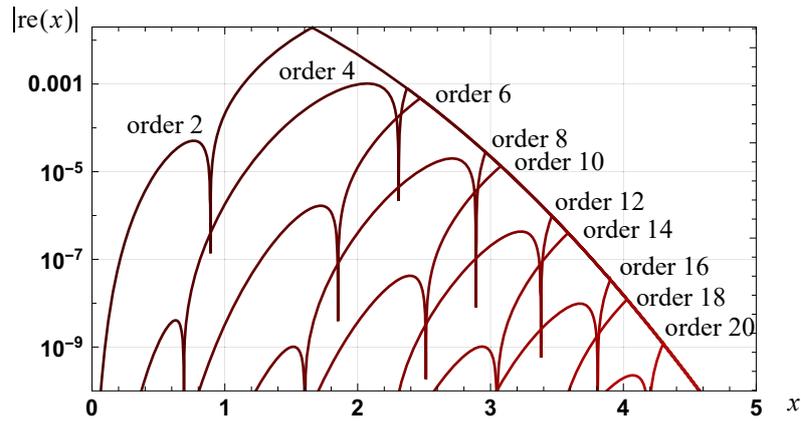

**Figure 7.** Graph of the relative errors in the, approximations $f_n$ to $\text{erf}(x)$, of orders 2, 4, 6, …, 20, based on utilizing the approximation $\text{erf}(x) \approx 1$ in an optimum manner.

where $T_n$ is a $n$th order Taylor series, $n$ odd, as specified in Table 1. The optimum transition points and the relative error bounds, for selected orders, are detailed in Table 4. The variation of the relative errors, with order, are shown in Figure 8. The change in the optimum transition point can be approximated according to

$$x_o(n) \approx 0.932 + 0.0560n - 0.0003503n^2, \qquad 3 \leq n \leq 61, \tag{60}$$

but, again, as small variations in $x_o(n)$ can lead to significant changes in the maximum relative error in the approximation for the error function, precise values for $x_o(n)$ are preferable. The clear superiority in the convergence of the spline based series is evident by a visual comparison of the relative errors shown in Figure 7 and Figure 8.

**Table 4.** The transition points, and the resulting relative error bounds, for Taylor series approximations specified by Equation 59. The transition points are based on sampling the interval $[0, 4]$ with 10000 points.

| Order: $n$ | Transition point: $x_o(n)$ | Relative error bound in $T_n$ |
|---|---|---|
| 1 | 0.8864 | 0.266 |
| 3 | 1.078 | 0.146 |
| 5 | 1.222 | 0.0917 |
| 7 | 1.344 | 0.0609 |
| 9 | 1.4532 | 0.0416 |
| 13 | 1.6452 | 0.0204 |
| 17 | 1.8144 | 0.0105 |
| 21 | 1.9672 | $5.44 \times 10^{-3}$ |
| 25 | 2.1084 | $2.89 \times 10^{-3}$ |
| 29 | 2.24 | $1.55 \times 10^{-3}$ |
| 37 | 2.4812 | $4.53 \times 10^{-4}$ |
| 45 | 2.70 | $1.35 \times 10^{-4}$ |
| 53 | 2.902 | $4.09 \times 10^{-5}$ |
| 61 | 3.09 | $1.24 \times 10^{-5}$ |

## 4   Variable Sub-interval Approximations for Error Function

An improved analytic approximation for the error function can be achieved by demarcating the interval $[0, x]$ into variable sub-intervals, e.g. the sub-intervals $[0, x/4]$, $[x/4, x/2]$, $[x/2, 3x/4]$ and $[3x/4, x]$ for the four sub-interval case, and by utilizing spline based integral approximations for each sub-interval. Chiani, 2002, utilized sub-intervals to enhance approximations for the complementary error function.

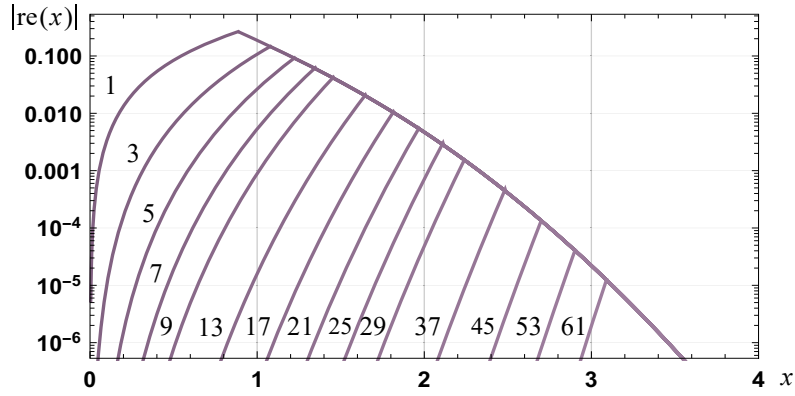

**Figure 8.** Graph of the magnitude of the relative error in Taylor series approximations to $\mathrm{erf}(x)$ that utilize an optimized change to the approximation $\mathrm{erf}(x) \approx 1$.

### Theorem 4.1   Variable Sub-Interval Approximations for Error Function

The $n$th order spline based approximation to the error function, based on $m$ equal width sub-intervals, is

$$f_{n,m}(x) = \frac{2}{\sqrt{\pi}} \sum_{i=0}^{m-1} \left[ \sum_{k=0}^{n} c_{n,k} \left[\frac{x}{m}\right]^{k+1} \left[ p\left[k, \frac{ix}{m}\right] \exp\left[-\frac{i^2 x^2}{m^2}\right] + (-1)^k p\left[k, \frac{(i+1)x}{m}\right] \exp\left[-\frac{(i+1)^2 x^2}{m^2}\right] \right] \right] \qquad (61)$$

where

$$p(k, x) = p^{(1)}(k-1, x) - 2x p(k-1, x), \qquad p(0, x) = 1. \qquad (62)$$

An alternative form is

$$f_{n,m}(x) = \frac{2}{\sqrt{\pi}} \left[ p_{n,0}(x) + \sum_{i=1}^{m-1} p_{n,i}(x) \exp\left[-\frac{i^2 x^2}{m^2}\right] + p_{n,m}(x) \exp(-x^2) \right] \qquad (63)$$

where

$$p_{n,0}(x) = \sum_{k=0}^{n} \frac{c_{n,k}}{m^{k+1}} \cdot p(k, 0) x^{k+1}, \qquad p_{n,m}(x) = \sum_{k=0}^{n} \frac{c_{n,k}}{m^{k+1}} \cdot (-1)^k p(k, x) x^{k+1},$$

$$p_{n,i}(x) = \sum_{k=0}^{n} \frac{c_{n,k}}{m^{k+1}} \cdot [1 + (-1)^k] p\left[k, \frac{ix}{m}\right] x^{k+1}. \qquad (64)$$

### Proof

The first result follows by applying Equation 27 in Theorem 2.1 to the sub-intervals $[0, x/m]$, $[x/m, 2x/m]$, …, $[x - x/m, x]$. The alternative form arises by expanding the outer summation in Equation 61 and collecting terms of similar form.

### 4.1   Explicit Expressions

A first order approximation, based on $m$ sub-intervals, is

$$f_{1,m}(x) = \frac{x}{m\sqrt{\pi}} \left[ 1 + 2 \sum_{i=1}^{m-1} \exp\left[\frac{-i^2 x^2}{m^2}\right] + \exp(-x^2) \right] + \frac{x^3}{3m^2 \sqrt{\pi}} \cdot \exp(-x^2). \qquad (65)$$

For the four sub-interval case, explicit expressions are

$$f_{0,4}(x) = \frac{x}{4\sqrt{\pi}}\left[1 + 2\exp\left[\frac{-x^2}{16}\right] + 2\exp\left[\frac{-x^2}{4}\right] + 2\exp\left[\frac{-9x^2}{16}\right] + \exp(-x^2)\right] \quad (66)$$

$$f_{1,4}(x) = \frac{x}{4\sqrt{\pi}}\left[1 + 2\exp\left[\frac{-x^2}{16}\right] + 2\exp\left[\frac{-x^2}{4}\right] + 2\exp\left[\frac{-9x^2}{16}\right] + \exp(-x^2)\right] + \frac{x^3}{48\sqrt{\pi}} \cdot \exp(-x^2) \quad (67)$$

Using the alternative form, a fourth order expression is

$$f_{4,4}(x) = \frac{x}{4\sqrt{\pi}}\left[1 - \frac{x^2}{288} + \frac{x^4}{322560}\right] +$$

$$\frac{x}{2\sqrt{\pi}} \cdot \exp\left[\frac{-x^2}{16}\right]\left[1 - \frac{x^2}{288} + \frac{47x^4}{107520} - \frac{x^6}{1,290,040} + \frac{x^8}{61,931,520}\right] +$$

$$\frac{x}{2\sqrt{\pi}} \cdot \exp\left[\frac{-x^2}{4}\right]\left[1 - \frac{x^2}{288} + \frac{187x^4}{107520} - \frac{x^6}{322560} + \frac{x^8}{3,870,720}\right] + \quad (68)$$

$$\frac{x}{2\sqrt{\pi}} \cdot \exp\left[\frac{-9x^2}{16}\right]\left[1 - \frac{x^2}{288} + \frac{1261x^4}{322560} - \frac{x^6}{143360} + \frac{3x^8}{2,293,760}\right] +$$

$$\frac{x}{4\sqrt{\pi}} \cdot \exp(-x^2)\left[1 + \frac{31x^2}{288} + \frac{101x^4}{15360} + \frac{19x^6}{80640} + \frac{x^8}{241920}\right]$$

A fourth order spline approximation, which utilizes sixteen sub-intervals, is detailed in Appendix 4. This expression, when utilized with the transition point $x_o = 7.1544$, yields an approximation with a relative error bound of $4.82 \times 10^{-16}$.

### 4.1.1 Results

The relative errors in the spline approximations of orders one to six, and for the case of four equal sub-intervals $[0, x/4]$, $[x/4, x/2]$, $[x/2, 3x/4]$ and $[3x/4, x]$, are shown in Figure 9.

### 4.2 Improved Approximation

The spline approximations utilizing variable sub-intervals can be improved by using the transition to the approximation $\text{erf}(x) \approx 1$ at a suitable point as specified by Equation 56. The relative error in the spline approximations of orders one to seven, and for the case of four equal sub-intervals $[0, x/4]$, $[x/4, x/2]$, $[x/2, 3x/4]$ and $[3x/4, x]$, are updated in Figure 10 to show the improvement associated with utilizing the optimum transition point to the approximation $\text{erf}(x) \approx 1$. The relative error bounds, and transition points, are detailed in Table 5 and Table 6 for the cases of four and sixteen sub-intervals.

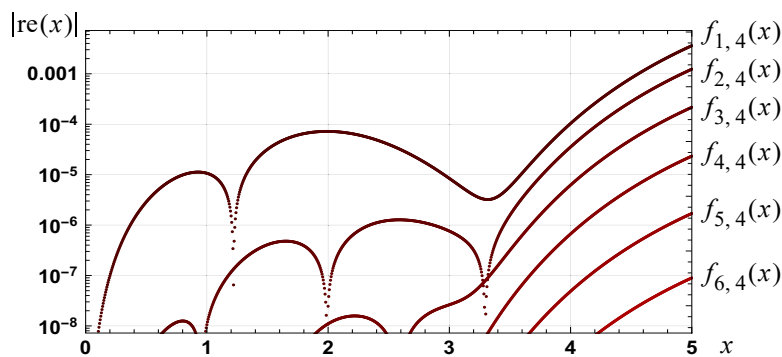

**Figure 9.** Graph of the relative errors in spline approximations to $\text{erf}(x)$, of orders one to six and based on four variable sub-intervals of equal width.

**Table 5.**   Transition point and relative error bound for the four equal sub-interval case. The transition points are based on sampling the interval [0, 8] with 10000 points.

| spline order | transition point | relative error bound |
|---|---|---|
| 0 | 2.7016 | $5.32 \times 10^{-3}$ |
| 1 | 3.292 | $7.21 \times 10^{-5}$ |
| 2 | 3.4544 | $1.27 \times 10^{-6}$ |
| 4 | 3.7208 | $1.43 \times 10^{-7}$ |
| 8 | 4.6616 | $4.34 \times 10^{-11}$ |
| 12 | 5.6784 | $9.75 \times 10^{-16}$ |
| 16 | 6.3736 | $2.01 \times 10^{-19}$ |
| 20 | 7.1544 | $4.62 \times 10^{-24}$ |
| 24 | 7.7136 | $1.06 \times 10^{-27}$ |

**Table 6.**   Transition point and relative error bound for the sixteen equal sub-interval case. The transition points are based on sampling the interval [0, 12] with 10000 points.

| spline order | transition point | relative error bound |
|---|---|---|
| 0 | 5.5008 | $3.32 \times 10^{-4}$ |
| 1 | 6.8796 | $2.82 \times 10^{-7}$ |
| 2 | 7.0224 | $3.14 \times 10^{-10}$ |
| 4 | 7.1544 | $4.82 \times 10^{-16}$ |
| 8 | 7.5996 | $6.22 \times 10^{-27}$ |
| 12 | 8.2032 | $4.16 \times 10^{-31}$ |
| 16 | 8.9244 | $1.66 \times 10^{-36}$ |
| 20 | 9.7284 | $4.68 \times 10^{-43}$ |
| 24 | 10.584 | $1.21 \times 10^{-50}$ |

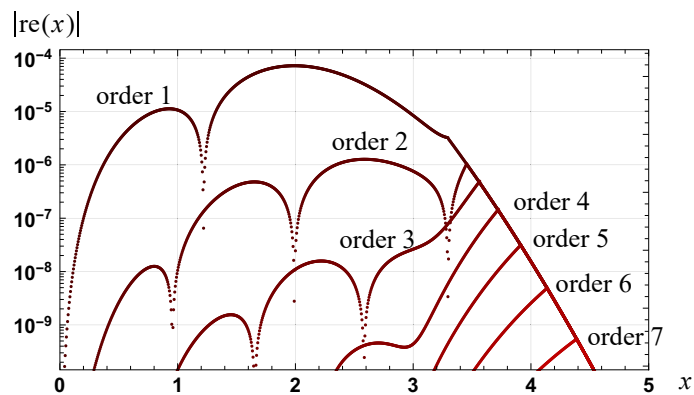

**Figure 10.**   Graph of the relative errors in approximations to $\text{erf}(x)$: first to seventh order spline based series based on four sub-intervals of equal width and with utilization of the approximation $\text{erf}(x) \approx 1$ at the optimum transition point.

### 4.2.1   Examples

A first order approximation, based on $m$ sub-intervals, as specified by Equation 65, yields the relative error bound of $7.21 \times 10^{-5}$ for four sub-intervals and with the transition point $x_o = 3.2928$; $4.51 \times 10^{-6}$ for eight sub-intervals and with the transition point $x_o = 4.784$; $2.82 \times 10^{-7}$ for sixteen sub-intervals and with the transition point of $x_o = 6.88$; and $1.10 \times 10^{-9}$ for sixty four sub-intervals and the transition point $x_o = 15.7888$.

The fourth order approximation, based on four equal sub-intervals, as specified by Equation 68, leads to the relative error bound of $1.43 \times 10^{-7}$ when used with the transition point $x_o = 3.7208$. A sixteenth order approximation, based on four equal sub-intervals, leads to an error bound of $2.01 \times 10^{-19}$ when used with the optimum transition point of $x_o = 6.3736$.

## 5   Dynamic Constant plus Spline Approximation

Consider the demarcation of the areas, as illustrated in Figure 11 and based on a resolution $\Delta$, that define the error function. It follows that

$$\text{erf}(x) = \sum_{k=0}^{\lfloor x/\Delta \rfloor} c_k + \frac{2}{\sqrt{\pi}} \int_{\Delta \lfloor x/\Delta \rfloor}^{x} e^{-\lambda^2} d\lambda \qquad \begin{cases} c_0 = 0 \\ c_k = \text{erf}(k\Delta) - \text{erf}[(k-1)\Delta], & k \geq 1. \end{cases} \tag{69}$$

For the general case of non-uniformly spaced intervals, as defined by the set of monotonically increasing points $\{x_0, x_1, x_2, \ldots, x_m\}$, and where it is not necessarily the case that $x > x_m$, the error function is defined according to

$$\text{erf}(x) = \sum_{k=1}^{m} c_k u(x - x_k) + \frac{2}{\sqrt{\pi}} \int_{x_S}^{x} e^{-\lambda^2} d\lambda \tag{70}$$

where $c_0 = 0$, $x_0 = 0$ and

$$c_k = \text{erf}(x_k) - \text{erf}[x_{k-1}], \qquad x_S = \sum_{k=1}^{m} [x_k - x_{k-1}] u(x - x_k). \tag{71}$$

A spline based approximation, as defined by Equation 27, can be utilized for the unknown integrals in Equation 69 and Equation 70. This leads to the following results:

### Theorem 5.1   Error Function Approximation: Dynamic Constant Plus a Spline Approximation

The error function, as defined by Equation 69 and Equation 70 can be approximated, respectively, by the approximations

$$f_{n,\Delta}(x) = \sum_{k=0}^{\lfloor x/\Delta \rfloor} c_k + \frac{2}{\sqrt{\pi}} \left[ \sum_{k=0}^{n} c_{n,k} \left[ x - \Delta \left\lfloor \frac{x}{\Delta} \right\rfloor \right]^{k+1} \left[ p\left[k, \Delta \left\lfloor \frac{x}{\Delta} \right\rfloor \right] \exp\left[-\Delta^2 \left\lfloor \frac{x}{\Delta} \right\rfloor^2\right] + (-1)^k p(k,x) \exp(-x^2) \right] \right] \tag{72}$$

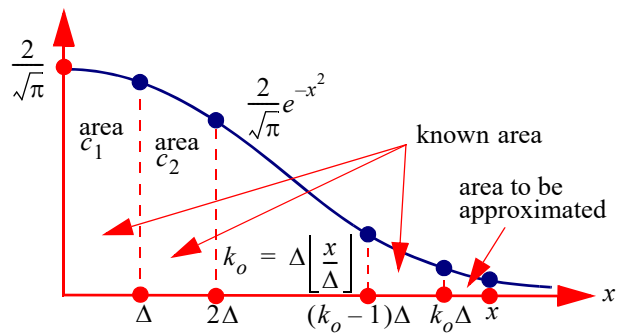

**Figure 11.**   Illustration of areas comprising $\text{erf}(x)$.

$$\text{erf}(x) \approx \sum_{k=1}^{m} c_k u(x - x_k) + \frac{2}{\sqrt{\pi}} \left[ \sum_{k=0}^{n} c_{n,k}(x - x_S)^{k+1} [p(k, x_S)\exp(-x_S^2) + (-1)^k p(k, x)\exp(-x^2)] \right] \quad (73)$$

**Proof**

These results arise from spline approximation of order $n$, as defined by Equation 27, for the integrals, respectively, over the intervals $[\Delta \lfloor x/\Delta \rfloor, x]$ and $[x_S, x]$.

### 5.1   Approximations of Orders Zeros to Four

Approximations of orders zero to four arising, from Theorem 5.1, are:

$$f_{0,\Delta}(x) = \sum_{k=0}^{\lfloor x/\Delta \rfloor} c_k + \frac{x - \Delta \lfloor x/\Delta \rfloor}{\sqrt{\pi}} \cdot \left[ e^{-\Delta^2 \lfloor x/\Delta \rfloor^2} + e^{-x^2} \right] \quad (74)$$

$$f_{1,\Delta}(x) = \sum_{k=0}^{\lfloor x/\Delta \rfloor} c_k + \frac{x - \Delta \lfloor x/\Delta \rfloor}{\sqrt{\pi}} \cdot \left[ e^{-\Delta^2 \lfloor x/\Delta \rfloor^2} + e^{-x^2} \right] - \frac{(x - \Delta \lfloor x/\Delta \rfloor)^2}{3\sqrt{\pi}} \cdot \left[ \Delta \lfloor x/\Delta \rfloor e^{-\Delta^2 \lfloor x/\Delta \rfloor^2} - xe^{-x^2} \right] \quad (75)$$

$$f_{2,\Delta}(x) = \sum_{k=0}^{\lfloor x/\Delta \rfloor} c_k + \frac{x - \Delta \lfloor x/\Delta \rfloor}{\sqrt{\pi}} \cdot \left[ e^{-\Delta^2 \lfloor x/\Delta \rfloor^2} + e^{-x^2} \right] - \frac{2(x - \Delta \lfloor x/\Delta \rfloor)^2}{5\sqrt{\pi}} \cdot \left[ \Delta \lfloor x/\Delta \rfloor e^{-\Delta^2 \lfloor x/\Delta \rfloor^2} - xe^{-x^2} \right] -$$
$$\frac{(x - \Delta \lfloor x/\Delta \rfloor)^3}{30\sqrt{\pi}} \cdot \left[ [1 - 2\Delta^2 \lfloor x/\Delta \rfloor^2] e^{-\Delta^2 \lfloor x/\Delta \rfloor^2} + [1 - 2x^2] e^{-x^2} \right] \quad (76)$$

$$f_{3,\Delta}(x) = \sum_{k=0}^{\lfloor x/\Delta \rfloor} c_k + \frac{x - \Delta \lfloor x/\Delta \rfloor}{\sqrt{\pi}} \cdot \left[ e^{-\Delta^2 \lfloor x/\Delta \rfloor^2} + e^{-x^2} \right] - \frac{3(x - \Delta \lfloor x/\Delta \rfloor)^2}{7\sqrt{\pi}} \cdot \left[ \Delta \lfloor x/\Delta \rfloor e^{-\Delta^2 \lfloor x/\Delta \rfloor^2} - xe^{-x^2} \right] -$$
$$\frac{(x - \Delta \lfloor x/\Delta \rfloor)^3}{21\sqrt{\pi}} \cdot \left[ [1 - 2\Delta^2 \lfloor x/\Delta \rfloor^2] e^{-\Delta^2 \lfloor x/\Delta \rfloor^2} + [1 - 2x^2] e^{-x^2} \right] + \quad (77)$$
$$\frac{(x - \Delta \lfloor x/\Delta \rfloor)^4}{70\sqrt{\pi}} \cdot \left[ \Delta \lfloor x/\Delta \rfloor \left[ 1 - \frac{2\Delta^2 \lfloor x/\Delta \rfloor^2}{3} \right] e^{-\Delta^2 \lfloor x/\Delta \rfloor^2} - x \left[ 1 - \frac{2x^2}{3} \right] e^{-x^2} \right]$$

$$f_{4,\Delta}(x) = \sum_{k=0}^{\lfloor x/\Delta \rfloor} c_k + \frac{x - \Delta \lfloor x/\Delta \rfloor}{\sqrt{\pi}} \cdot \left[ e^{-\Delta^2 \lfloor x/\Delta \rfloor^2} + e^{-x^2} \right] - \frac{4(x - \Delta \lfloor x/\Delta \rfloor)^2}{9\sqrt{\pi}} \cdot \left[ \Delta \lfloor x/\Delta \rfloor e^{-\Delta^2 \lfloor x/\Delta \rfloor^2} - xe^{-x^2} \right] -$$
$$\frac{(x - \Delta \lfloor x/\Delta \rfloor)^3}{18\sqrt{\pi}} \cdot \left[ [1 - 2\Delta^2 \lfloor x/\Delta \rfloor^2] e^{-\Delta^2 \lfloor x/\Delta \rfloor^2} + [1 - 2x^2] e^{-x^2} \right] +$$
$$\frac{(x - \Delta \lfloor x/\Delta \rfloor)^4}{42\sqrt{\pi}} \cdot \left[ \Delta \lfloor x/\Delta \rfloor \left[ 1 - \frac{2\Delta^2 \lfloor x/\Delta \rfloor^2}{3} \right] e^{-\Delta^2 \lfloor x/\Delta \rfloor^2} - x \left[ 1 - \frac{2x^2}{3} \right] e^{-x^2} \right] + \quad (78)$$
$$\frac{(x - \Delta \lfloor x/\Delta \rfloor)^5}{1260\sqrt{\pi}} \cdot \left[ \left[ 1 - 4\Delta^2 \lfloor x/\Delta \rfloor^2 + \frac{4\Delta^4 \lfloor x/\Delta \rfloor^4}{3} \right] e^{-\Delta^2 \lfloor x/\Delta \rfloor^2} + \left[ 1 - 4x^2 + \frac{4x^4}{3} \right] e^{-x^2} \right]$$

## 5.2 Results

For a resolution of $\Delta = 1/2$, the coefficients are tabulated in Table 7.

**Table 7.** Coefficient values for the case of $\Delta = 1/2$.

| $k$ | Definition for $c_k$ | $c_k$ |
|---|---|---|
| 1 | $\mathrm{erf}(1/2)$ | $5.204998778 \times 10^{-1}$ |
| 2 | $\mathrm{erf}(1) - \mathrm{erf}(1/2)$ | $3.222009151 \times 10^{-1}$ |
| 3 | $\mathrm{erf}(3/2) - \mathrm{erf}(1)$ | $1.234043535 \times 10^{-1}$ |
| 4 | $\mathrm{erf}(2) - \mathrm{erf}(3/2)$ | $2.921711854 \times 10^{-2}$ |
| 5 | $\mathrm{erf}(5/2) - \mathrm{erf}(2)$ | $4.270782964 \times 10^{-3}$ |
| 6 | $\mathrm{erf}(3) - \mathrm{erf}(5/2)$ | $3.848615204 \times 10^{-4}$ |
| 7 | $\mathrm{erf}(7/2) - \mathrm{erf}(3)$ | $2.134739863 \times 10^{-5}$ |
| 8 | $\mathrm{erf}(4) - \mathrm{erf}(7/2)$ | $7.276811144 \times 10^{-7}$ |
| 9 | $\mathrm{erf}(9/2) - \mathrm{erf}(4)$ | $1.522064186 \times 10^{-8}$ |
| 10 | $\mathrm{erf}(5) - \mathrm{erf}(9/2)$ | $1.950785844 \times 10^{-10}$ |
| 11 | $\mathrm{erf}(11/2) - \mathrm{erf}(5)$ | $1.530101947 \times 10^{-12}$ |
| 12 | $\mathrm{erf}(6) - \mathrm{erf}(11/2)$ | $7.336328181 \times 10^{-15}$ |

A resolution of $\Delta = 1/2$ yields a relative error bound of $1.16 \times 10^{-5}$ for a second order approximation, a relative error bound of $1.35 \times 10^{-9}$ for a fourth order approximation, a relative error bound of $7.15 \times 10^{-14}$ for a sixth order approximation and a relative error bound of $9.03 \times 10^{-37}$ for a sixteenth order approximation. These bounds are based on 10000 equal spaced samples in the interval $[0, 8]$.

The variation of the relative error bound with resolution, and order, is detailed in Figure 12. The nature of the variation of the relative error, for orders two, three and four, is shown in Figure 13 for the case of resolution of 0.5. It is possible to obtain better results by using non-uniformly spaced intervals but the improvement, in general, does not warrant the increase in complexity.

## 6 A Dynamical System to Yield Improved Approximations

It is possible to utilize the approximations detailed in Theorem 2.1 and Theorem 4.1 as the basis for determining new approximations with a lower relative error. The approach is indirect and based on considering the feedback system illustrated in Figure 14 which has dynamically varying feedback. The differential equation characterizing the system is

$$y'(t) + f_M(t)y(t) = x(t). \tag{79}$$

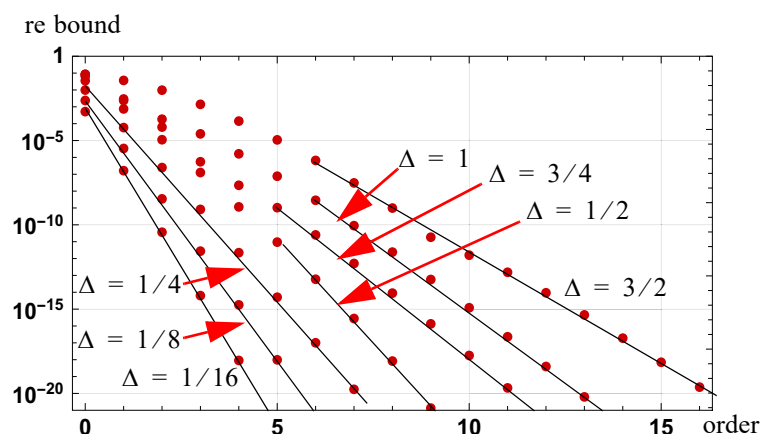

**Figure 12.** Graph of the relative error bound, versus the order of approximation, for various set resolutions.

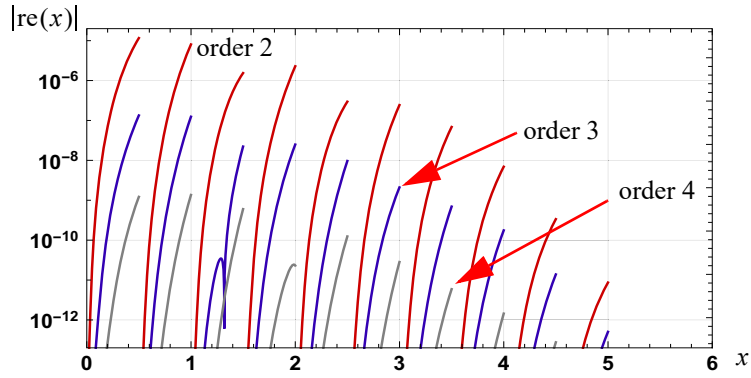

**Figure 13.**   Graph of the relative errors, based on a resolution of $\Delta = 0.5$, in second to fourth order approximations to $\text{erf}(x)$.

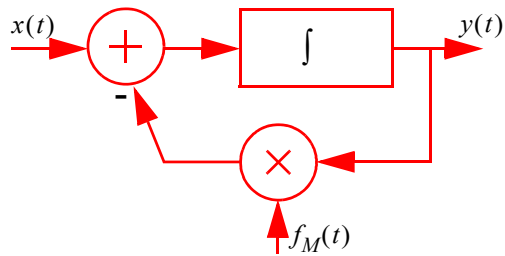

**Figure 14.**   Feedback system with dynamically varying (modulated) feedback.

For specific input, $x$, and modulated feedback, $f_M$, signals the output has a known form. For example, for the case of $x(t) = f_M(t) = \text{erf}(t)u(t)$, the output signal, assuming zero initial conditions, is

$$y(t) = 1 - \exp\left[\frac{1}{\sqrt{\pi}} \cdot \left[1 - e^{-t^2}\right] - t\,\text{erf}(t)\right], \qquad t \geq 0. \tag{80}$$

For the case of

$$x(t) = f_M(t) = \frac{4}{\sqrt{\pi}} e^{-t^2} \text{erf}(t) u(t) \tag{81}$$

the output signal, assuming zero initial conditions, is

$$y(t) = 1 - \exp[-\text{erf}^2(t)], \qquad t \geq 0. \tag{82}$$

This case facilitates approximations for the error function, which can be made arbitrarily accurate and which are valid for the positive real line.

### Theorem 6.1   Dynamical System Approximations for Error Function

Based on the differential equation specified by Equation 79, a $n$th order approximation to $\text{erf}(x)$, for the case of $x \geq 0$, can be defined according to

$$f_n(x) = \sqrt{p_{n,0} + p_{n,1}(x)e^{-x^2} + p_{n,2}(x)e^{-2x^2}} \tag{83}$$

where, for the case of $n \geq 2$:

$$p_{n,1}(x) = \alpha_0 + \alpha_2 x^2 + \ldots + \alpha_m x^m, \qquad m = \begin{cases} n-1 & n \text{ odd} \\ n & m \text{ even} \end{cases}$$

$$p_{n,2}(x) = \beta_0 + \beta_2 x^2 + \ldots + \beta_{2n} x^{2n}$$

$$p_{n,0} = -(\alpha_0 + \beta_0)$$

(84)

and with

$$\alpha_m = \frac{-4}{\pi} \cdot c_{n,m} a_{m,0}$$

$$\alpha_{m-2i} = \frac{(m-2i+2)\alpha_{m-2i+2}}{2} - \frac{4}{\pi} \cdot c_{n,m-2i} a_{m-2i,0}, \qquad i \in \left\{1, \ldots, \frac{m}{2} - 1\right\}$$

(85)

$$\alpha_0 = \alpha_2 - \frac{4}{\pi} \cdot c_{n,0} a_{0,0}$$

$$\beta_{2n} = \frac{-2}{\pi} \cdot (-1)^n c_{n,n} a_{n,n}$$

$$\beta_{2n-2i} = \frac{(n-i+1)\beta_{2n-2i+2}}{2} - \frac{2}{\pi} \sum_{k=n-i}^{\min\{2n-2i,\,n\}} (-1)^k c_{n,k} a_{k,2(n-i)-k}, \qquad i \in \{1, \ldots, n-1\}$$

(86)

$$\beta_0 = \frac{\beta_2}{2} - \frac{2}{\pi} \cdot c_{n,0} a_{0,0}$$

Here, the coefficients $a_{i,j}$, $i, j \in \{0, 1, \ldots, n\}$ are defined by the expansion

$$p(k, x) = a_{k,0} + a_{k,1} x + a_{k,2} x^2 + \ldots + a_{k,k} x^k, \qquad k \in \{0, 1, \ldots, n\},$$

(87)

arising from the polynomials (Equation 26)

$$p(k, x) = p^{(1)}(k-1, x) - 2xp[k-1, x], \qquad p(0, x) = 1.$$

(88)

Finally, it is the case that

$$\lim_{n \to \infty} f_n(x) = \text{erf}(x), \qquad x \geq 0,$$

(89)

with the convergence being uniform.

**Proof**

The proof is detailed in Appendix 5.

### 6.1    Explicit Approximations

Explicit approximations for orders zero to four for $\text{erf}(x)$, $x \geq 0$ are:

$$f_0(x) = \frac{1}{\sqrt{\pi}} \sqrt{3 - 2e^{-x^2} - e^{-2x^2}}$$

(90)

$$f_1(x) = \frac{1}{\sqrt{\pi}} \sqrt{\frac{19}{6} - 2e^{-x^2} - \frac{7e^{-2x^2}}{6}\left[1 + \frac{2x^2}{7}\right]}$$

(91)

$$f_2(x) = \frac{1}{\sqrt{\pi}} \sqrt{\frac{63}{20} - \frac{29e^{-x^2}}{15}\left[1 - \frac{x^2}{29}\right] - \frac{73e^{-2x^2}}{60}\left[1 + \frac{26x^2}{73} + \frac{4x^4}{73}\right]}$$

(92)

$$f_3(x) = \frac{1}{\sqrt{\pi}} \sqrt{\frac{22}{7} - \frac{40e^{-x^2}}{21}\left[1 - \frac{x^2}{20}\right] - \frac{26e^{-2x^2}}{21}\left[1 + \frac{10x^2}{26} + \frac{x^4}{13} + \frac{x^6}{130}\right]}$$

(93)

$$f_4(x) = \frac{1}{\sqrt{\pi}}\sqrt{\frac{377}{120} - \frac{596e^{-x^2}}{315}\left[1 - \frac{17x^2}{298} + \frac{x^4}{1192}\right] - \frac{3149e^{-2x^2}}{2520}\left[1 + \frac{1258x^2}{3149} + \frac{278x^4}{3149} + \frac{112x^6}{9447} + \frac{8x^8}{9447}\right]} \quad (94)$$

## 6.2   Results

The relative error bounds associated with the approximations to erf(x), are detailed in Table 8. The graphs of the relative errors in the approximations are shown in Figure 15. The clear advantage of the proposed approximations is evident with the improvement increasing with the order of the initial approximation (i.e. a function with an initial lower relative error bound leads to an increasingly lower relative error bound). The other clear advantage of the approximations, as is evident in Figure 15, is that the relative error is bounded as $x \to \infty$.

**Table 8.** Relative error bounds, over the interval $(0, \infty)$, for approximations to erf(x) as defined in Theorem 6.1.

| Order of approx. | Relative error bound: original series - optimum transition point (Table 3) | Relative error bound: approx. defined by Equation 83 |
|---|---|---|
| 0 | 0.0851 | $2.68 \times 10^{-2}$ |
| 1 | 0.0362 | $3.98 \times 10^{-3}$ |
| 2 | 0.0195 | $1.34 \times 10^{-3}$ |
| 3 | $7.36 \times 10^{-3}$ | $2.03 \times 10^{-4}$ |
| 4 | $1.03 \times 10^{-3}$ | $1.82 \times 10^{-5}$ |
| 6 | $4.75 \times 10^{-4}$ | $9.20 \times 10^{-7}$ |
| 8 | $2.79 \times 10^{-5}$ | $1.69 \times 10^{-8}$ |
| 10 | $1.35 \times 10^{-5}$ | $7.43 \times 10^{-10}$ |
| 12 | $9.78 \times 10^{-7}$ | $1.67 \times 10^{-11}$ |
| 14 | $4.00 \times 10^{-7}$ | $6.47 \times 10^{-13}$ |
| 16 | $3.44 \times 10^{-8}$ | $1.68 \times 10^{-14}$ |
| 18 | $1.22 \times 10^{-8}$ | $5.90 \times 10^{-16}$ |
| 20 | $1.20 \times 10^{-9}$ | $1.73 \times 10^{-17}$ |
| 22 | $3.76 \times 10^{-10}$ | $5.56 \times 10^{-19}$ |
| 24 | $4.18 \times 10^{-11}$ | $1.79 \times 10^{-20}$ |

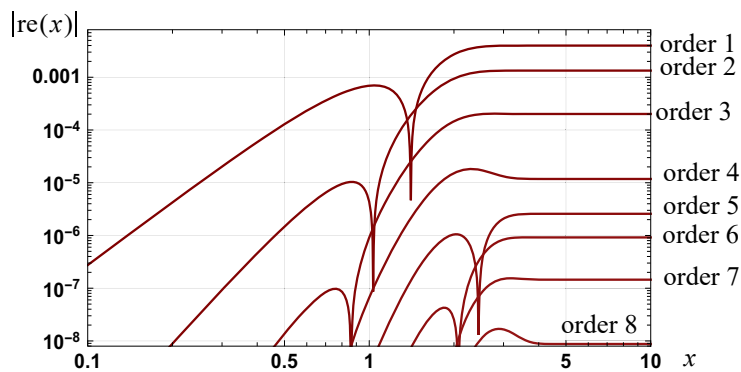

**Figure 15.** Graph of the relative errors in approximations, of orders one to eight, to erf(x) as defined in Theorem 6.1.

### 6.3    Extension

By utilizing the approximations detailed in Theorem 4.1 similar approximations can be detailed, with lower relative error. For example, the first order approximation, $f_{1,4}$, which is based on four equal sub-intervals and is defined by Equation 67, yields the approximation

$$f_{1,4}(x) = \frac{1}{\sqrt{\pi}}\sqrt{\frac{128177}{40800} - \frac{e^{-x^2}}{2} - \frac{16e^{-17x^2/16}}{17} - \frac{4e^{-5x^2/4}}{5} - \frac{16e^{-25x^2/16}}{25} - \frac{25e^{-2x^2}}{96}\left[1 + \frac{2x^2}{25}\right]} \quad (95)$$

which has a relative error bound of $2.83 \times 10^{-6}$. With an optimum transition point of 3.292 the original approximation has a relative error bound of $7.21 \times 10^{-5}$ (see Table 5).

### 6.4    Notes

First, the constants $p_{n,0}$, $n \in \{0, 1, \ldots\}$, as defined in Equation 84, form a series that in the limit converges to 1. It then follows that the corresponding series converges to $\pi$:

$$3, \frac{19}{6}, \frac{63}{20}, \frac{22}{7}, \frac{377}{120}, \frac{174169}{55440}, \frac{4,528,409}{1,441,440}, \ldots \quad (96)$$

Second, the square root functional structure has been utilized for approximations to the error function as is evident from the approximations detailed in Table 1. It is easy to conclude that the form

$$f_n(x) = \sqrt{p_{n,0} - p_{n,1}(x)e^{-k_1 x^2} - p_{n,2}(x)e^{-k_2 x^2} - \ldots} \quad (97)$$

is well suited for approximating the error function.

## 7    Applications

This section details indicative applications of the approximations for the error function that have been detailed above.

The distinct analytical forms, that are specified in Theorem 2.1, Theorem 2.3, Theorem 4.1 and Theorem 6.1, for approximations to the error function, in general, facilitate analysis for different applications. For example, the form detailed in Theorem 2.1 underpins approximations to $\exp(-x^2)$ as detailed in Section 7.2. The form detailed in Theorem 6.1 underpins analytical approximations for the power associated with the output of a non-linearity modelled by the error function when subject to a sinusoidal signal. The approximations are detailed in Section 7.6.

For applications, where a set relative error bound over a set interval is required, the approximation that is appropriate will depend, in part, on the domain over which an approximation is required as well as the level of the relative error bound that is acceptable. For example, the approximations detailed in Theorem 2.1 and Theorem 2.3 lead to simple analytical forms and with modest relative error bounds over $[0, \infty)$ when used with an appropriate transition point to the approximation of $\text{erf}(x) \approx 1$. Without the use of a transition point, such approximations are likely to be best suited for a restricted domain, for example, the domain $[0, 3/\sqrt{2}]$ which is consistent with the three sigma case arising from a Gaussian distribution. The fourth order approximations, as specified by Equation 36 and Equation 54, have relative error bounds, respectively, of $1.02 \times 10^{-3}$ and $1.90 \times 10^{-4}$ over $[0, 3/\sqrt{2}]$. Section 7.1 provides examples of approximations that are consistent with set relative error bounds, over the interval $[0, \infty)$, of $10^{-4}$, $10^{-6}$, $10^{-10}$ and $10^{-16}$.

### 7.1    Error Function Approximations: Set Relative Error Bounds

Consider the case where an approximation for the error function, with a relative error bound over the positive real line of $10^{-4}$, is required. A 47th order Taylor series, with a transition point of $x_o = 2.752$, yields a relative bound of $1.00 \times 10^4$.

An eighth order spline approximation, with a transition point of $x_o = 2.963$, yields a relative error bound of $2.79 \times 10^{-5}$. The approximation, according to Equation 56, is

$$f_8(x) = \frac{x}{\sqrt{\pi}} \cdot u[x_o - x] \cdot \left[ 1 - \frac{7x^2}{102} + \frac{x^4}{340} - \frac{x^6}{18564} + \frac{x^8}{5,250,960} + \right.$$

$$\left[ 1 + \frac{41x^2}{102} + \frac{101x^4}{1020} + \frac{1591x^6}{92820} + \frac{4793x^8}{2,162,160} + \frac{2017x^{10}}{9,189,180} + \frac{38x^{12}}{2,297,295} + \right.$$
$$\left. \left. \frac{31x^{14}}{34,459,425} + \frac{x^{16}}{34,459,425} \right] e^{-x^2} \right] + \quad (98)$$

$$[1 - u[x_o - x]]$$

A seventh order approximation, with a transition point of $x_o = 2.65$, yields a relative error bound of $1.79 \times 10^{-4}$.

A first order spline approximation, based on four equal sub-intervals $[0, x/4]$, $[x/4, x/2]$, $[x/2, 3x/4]$ and $[3x/4, x]$, is defined according to

$$f_{1,4}(x) = \left[ \frac{x}{4\sqrt{\pi}} \cdot \left[ 1 + 2\exp\left[\frac{-x^2}{16}\right] + 2\exp\left[\frac{-x^2}{4}\right] + 2\exp\left[\frac{-9x^2}{16}\right] + \exp[-x^2] \right] + \frac{x^3}{48\sqrt{\pi}} \cdot \exp[-x^2] \right] u[x_o - x] + \quad (99)$$

$$[1 - u[x_o - x]]$$

and yields a relative error bound of $7.21 \times 10^{-5}$ with the transition point $x_o = 3.292$.

A dynamic constant plus a spline approximation of order 2, and based on a resolution of $\Delta = 19/20$ achieves a relative error bound of $8.33 \times 10^{-5}$ (10000 points in the interval $[0, 5]$). The approximation is

$$f_{2,\Delta}(x) = \sum_{k=0}^{\lfloor 20x/19 \rfloor} c_k + \frac{1}{\sqrt{\pi}}\left[x - \frac{19}{20} \cdot \left\lfloor \frac{20x}{19} \right\rfloor\right]\left[\exp\left[\frac{-361}{400}\left\lfloor \frac{20x}{19} \right\rfloor^2\right] + e^{-x^2}\right] -$$
$$\frac{2}{5\sqrt{\pi}}\left[x - \frac{19}{20} \cdot \left\lfloor \frac{20x}{19} \right\rfloor\right]^2\left[\frac{19}{20} \cdot \left\lfloor \frac{20x}{19} \right\rfloor \exp\left[\frac{-361}{400}\left\lfloor \frac{20x}{19} \right\rfloor^2\right] - xe^{-x^2}\right] - \quad (100)$$
$$\frac{1}{30\sqrt{\pi}}\left[x - \frac{19}{20} \cdot \left\lfloor \frac{20x}{19} \right\rfloor\right]^3 \cdot \left[\left[1 - \frac{361}{200}\left\lfloor \frac{20x}{19} \right\rfloor^2\right]\exp\left[\frac{-361}{400}\left\lfloor \frac{20x}{19} \right\rfloor^2\right] + (1 - 2x^2)e^{-x^2}\right]$$

where

$$c_0 = 0, \quad c_k = \text{erf}\left[\frac{19k}{20}\right] - \text{erf}\left[\frac{19(k-1)}{20}\right], \quad k \in \{1, 2, \ldots\},$$
$$c_1 = 0.82089081, \quad c_2 = 0.17189962, \quad (101)$$
$$c_3 = 7.1539145 \times 10^{-3}, \quad c_4 = 5.5579 \times 10^{-5}, \quad c_5 = c_6 = \ldots = 0.$$

Here, the approximation of $\text{erf}(x) \approx 1$, for $x \geq 57/20$ (after three intervals) can be utilized without impacting the relative error bound.

Utilizing a fourth order spline approximation and iteration consistent with Theorem 6.1, the approximation

$$f_4(x) = \frac{1}{\sqrt{\pi}}\sqrt{\frac{377}{120} - \frac{596e^{-x^2}}{315}\left[1 - \frac{17x^2}{298} + \frac{x^4}{1192}\right] - \frac{3149e^{-2x^2}}{2520}\left[1 + \frac{1258x^2}{3149} + \frac{278x^4}{3149} + \frac{112x^6}{9447} + \frac{8x^6}{9447}\right]} \quad (102)$$

yields a relative error bound of $1.82 \times 10^{-5}$.

Details of approximations that are consistent with higher order relative error bounds are detailed in Table 9.

**Table 9.** Approximations that are consistent with a set relative error bound. The actual relative error bound is specified by $re_B$.

| Relative error bound | Spline approx: Theorem 3.1 | Variable interval approx: Theorem 4.1 | Dynamic constant plus spline approx: Theorem 5.1 | Iterative approx: Theorem 6.1 |
|---|---|---|---|---|
| $10^{-6}$ | order $= 12$ <br> $x_o = 3.4625$ <br> $re_B = 9.78 \times 10^{-7}$ | order $= 5$ <br> 3 sub-intervals <br> $x_o = 3.51$ <br> $re_B = 6.96 \times 10^{-7}$ | order $= 3$ <br> resolution $= 3/4$ <br> $re_B = 5.53 \times 10^{-7}$ | order $= 6$ <br> $re_B = 9.20 \times 10^{-7}$ |
| $10^{-10}$ | order $= 23$ <br> $x_o = 4.581$ <br> $re_B = 9.31 \times 10^{-11}$ | order $= 8$ <br> 4 sub-intervals <br> $x_o = 4.6616$ <br> $re_B = 4.34 \times 10^{-11}$ | order $= 4$ <br> resolution $= 3/8$ <br> $re_B = 9.12 \times 10^{-11}$ | order $= 11$ <br> $re_B = 1.34 \times 10^{-10}$ <br> order $= 12$ <br> $re_B = 1.67 \times 10^{-11}$ |
| $10^{-16}$ | order $= 39$ <br> $x_o = 5.9017$ <br> $re_B = 7.21 \times 10^{-17}$ | order $= 11$ <br> 6 sub-intervals <br> $x_o = 5.98$ <br> $re_B = 2.75 \times 10^{-17}$ | order $= 6$ <br> resolution $= 1/4$ <br> $re_B = 1.01 \times 10^{-17}$ | order $= 19$ <br> $re_B = 1.18 \times 10^{-16}$ <br> order $= 20$ <br> $re_B = 1.73 \times 10^{-17}$ |

## 7.2 Approximation for Exp(-x²)

A $n$th order approximation to the Gaussian function $\exp(-x^2)$ is detailed in the following theorem:

### Theorem 7.1 Approximation for Gaussian Function

A $n$th order approximation, $g_n$ to the Gaussian function $\exp(-x^2)$ is

$$g_n(x) = \frac{\sum_{k=0}^{n} c_{n,k}(k+1)x^k p(k,0)}{1 + \sum_{k=0}^{n} c_{n,k}(-1)^{k+1}x^k[p(k,x)[k+1-2x^2)] + xp^{(1)}(k,x)]} \quad (103)$$

where $c_{n,k}$ is defined by Equation 21 and $p(k,x)$ is defined by Equation 26.

**Proof**

The proof is detailed in Appendix 6.

### 7.2.1 Approximations

Approximations to $\exp(-x^2)$, of orders zero to five, are:

$$g_0(x) = \frac{1}{1+2x^2} \qquad g_1(x) = \frac{1}{1+x^2+\frac{2x^4}{3}} \qquad g_2(x) = \frac{1-x^2/10}{1+\frac{9x^2}{10}+\frac{2x^4}{5}+\frac{2x^6}{15}} \quad (104)$$

$$g_3(x) = \frac{1 - x^2/7}{1 + \frac{6x^2}{7} + \frac{5x^4}{14} + \frac{2x^6}{21} + \frac{2x^8}{105}} \qquad g_4(x) = \frac{1 - x^2/6 + x^4/252}{1 + \frac{5x^2}{6} + \frac{85x^4}{252} + \frac{11x^6}{126} + \frac{x^8}{63} + \frac{2x^{10}}{945}} \qquad (105)$$

$$g_5(x) = \frac{1 - \frac{2x^2}{11} + \frac{x^4}{132}}{1 + \frac{9x^2}{11} + \frac{43x^4}{132} + \frac{x^6}{12} + \frac{x^8}{66} + \frac{x^{10}}{495} + \frac{2x^{12}}{10395}} \qquad (106)$$

### 7.2.2  Results

The relative errors in the above defined approximations to $\exp(-x^2)$ are detailed in Figure 16 for approximations of order 0, 2, 4, 6, 8, 10 and 12 along with the relative error in Taylor series for orders 1, 3, 5, …, 15. The clear superiority of the defined approximations is evident.

### 7.2.3  Comparison

The following $n$th order approximation for $\exp(-x^2)$ has been proposed in Howard 2019 (eqn. 77):

$$h_n(x) = \frac{1 - c_{n,0}x^2 + c_{n,1}x^4 - c_{n,2}x^6 + \ldots + c_{n,n}(-1)^{n+1}x^{2n+2}}{1 + c_{n,0}x^2 + c_{n,1}x^4 + c_{n,2}x^6 + \ldots + c_{n,n}x^{2n+2}} \qquad (107)$$

where $c_{n,k}$ is defined by Equation 21. The relative error bounds over the interval $[0, 3/\sqrt{2}]$ (the three sigma bound case for Gaussian probability distributions) for this approximation, and the approximation defined by Equation 103, are detailed in Table 10. The tabulated results clearly show that this approximation is more accurate the approximation detailed in Equation 103. The improvement is consistent with the higher order Padé approximant being used.

The following approximations (seventh and fifth order) yield relative error bounds of less than $0.001$ over the interval $[0, 3/\sqrt{2}]$:

$$g_7(x) = \frac{1 - \frac{x^2}{5} + \frac{x^4}{78} - \frac{x^6}{4290}}{1 + \frac{4x^2}{5} + \frac{61x^4}{195} + \frac{34x^6}{429} + \frac{83x^8}{5720} + \frac{x^{10}}{495} + \frac{7x^{12}}{32175} + \frac{4x^{14}}{225225} + \frac{2x^{16}}{2027025}} \qquad (108)$$

$$h_5(x) = \frac{1 - \frac{x^2}{2} + \frac{5x^4}{44} - \frac{x^6}{66} + \frac{x^8}{792} - \frac{x^{10}}{15840} + \frac{x^{12}}{665280}}{1 + \frac{x^2}{2} + \frac{5x^4}{44} + \frac{x^6}{66} + \frac{x^8}{792} + \frac{x^{10}}{15840} + \frac{x^{12}}{665280}} \qquad (109)$$

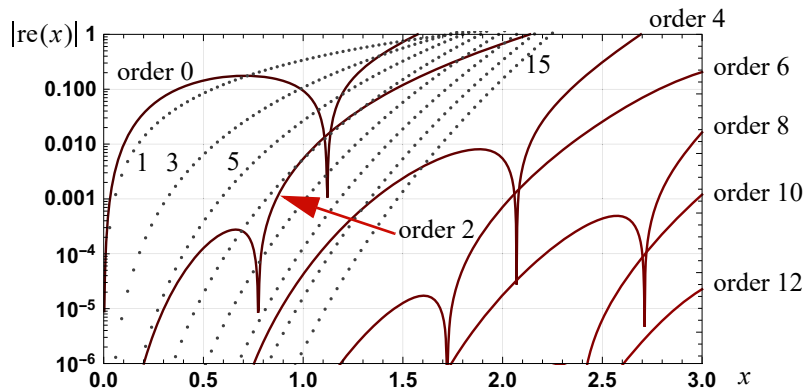

**Figure 16.** Graph of the magnitude of the relative errors in approximations to $\exp(-x^2)$, as defined by Equation 103, of orders 0, 2, 4, 6, 8, 10 and 12. The dotted curves are the relative errors associated with Taylor series of orders 1, 3, 5, 7, 9, 11, 13 and 15.

A twenty seventh order Taylor series approximation yields a relative bound of $1.03 \times 10^{-3}$.

Table 10. Relative error bounds for approximations over the interval $[0, 3/\sqrt{2}]$.

| Order of approx. | Relative error bound: Equation 103 | Relative error bound: Howard 2019, eqn. 77 |
|---|---|---|
| 0 | 8.00 | 35.6 |
| 1 | 3.74 | 6.98 |
| 2 | 0.957 | 0.767 |
| 3 | $1.25 \times 10^{-1}$ | $5.25 \times 10^{-2}$ |
| 4 | $8.04 \times 10^{-3}$ | $2.42 \times 10^{-3}$ |
| 5 | $7.71 \times 10^{-3}$ | $7.98 \times 10^{-5}$ |
| 6 | $2.09 \times 10^{-3}$ | $1.97 \times 10^{-6}$ |
| 7 | $3.72 \times 10^{-4}$ | $3.75 \times 10^{-8}$ |
| 8 | $5.02 \times 10^{-5}$ | $5.69 \times 10^{-10}$ |
| 10 | $4.54 \times 10^{-7}$ | $7.22 \times 10^{-14}$ |
| 12 | $1.09 \times 10^{-9}$ | $4.62 \times 10^{-18}$ |

## 7.3  Upper and Lower Bounded Approximations to Error Function

Establishing bounds for $\text{erf}(x)$ has received modest research interest, e.g. Alzer 2010 and published bounds for $\text{erf}(x)$ for the case of $x > 0$ include: Chu 1955:

$$\sqrt{1 - \exp[-px^2]} \leq \text{erf}(x) \leq \sqrt{1 - \exp[-qx^2]}, \qquad p \in (0, 1], q \in [4/\pi, \infty]. \tag{110}$$

corollary 4.2 of Neuman 2013:

$$\frac{2x}{\sqrt{\pi}} \exp\left[\frac{-x^2}{3}\right] \leq \text{erf}(x) \leq \frac{4x}{3\sqrt{\pi}}\left[1 + \frac{\exp(-x^2)}{2}\right] \tag{111}$$

and refinements to the form proposed by Chu, e.g. Yang 2016, 2018 (Corollary 3.4):

$$\sqrt{1 - \frac{20}{3\pi}\left[1 - \frac{\pi}{4}\right]\exp\left[\frac{-8x^2}{5}\right] - \frac{8}{3}\left[1 - \frac{5}{2\pi}\right]\exp(-x^2)} \leq \text{erf}(x) \leq$$

$$\sqrt{1 - \lambda(p_0)\exp(-p_0 x^2) - [1 - \lambda(p_0)]\exp[-\mu(p_0)x^2]} \tag{112}$$

where

$$p_0 = \frac{21\pi - 60 + \sqrt{3(147\pi^2 - 920\pi + 1440)}}{30(\pi - 3)}$$

$$\lambda(p) = \frac{4[7\pi - 20 - 5(\pi - 3)p]}{\pi(15p^2 - 40p + 28)}, \qquad \mu(p) = \frac{4(5p - 7)}{5(3p - 4)}. \tag{113}$$

The relative error in these bounds are detailed in Figure 17.

Utilizing the results of Lemma 1, it follows that any of the approximations detailed in Theorem 3.1, Theorem 4.1, Theorem 5.1 or Theorem 6.1 can be utilized to create upper and lower bounded functions for $\text{erf}(x)$, $x > 0$, of arbitrary accuracy and with an arbitrary relative error bound. For example, the approximation $f_{1,4}$ specified by Equation 99, yields the functional bounds:

$$\frac{f_{1,4}(x)}{1 + \varepsilon_B} \leq \text{erf}(x) \leq \frac{f_{1,4}(x)}{1 - \varepsilon_B}, \qquad \varepsilon_B = 7.21 \times 10^{-5}, x > 0, \tag{114}$$

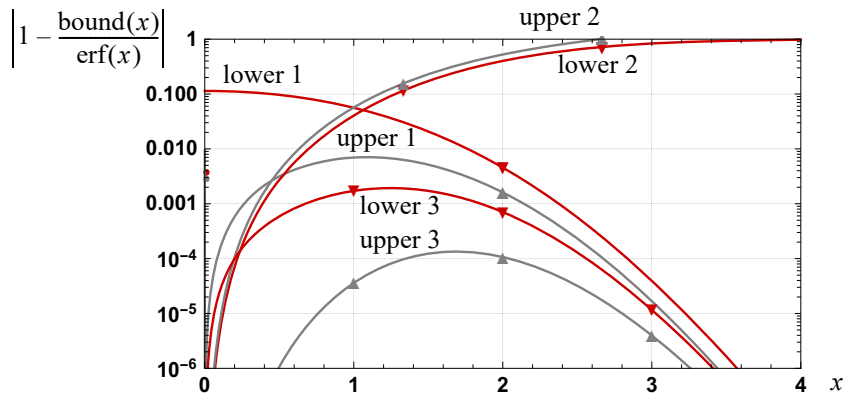

**Figure 17.** Relative error in upper and lower bounds to erf($x$) as, respectively, defined by Equation 110, Equation 111 and Equation 112. The parameters $p = 1$ and $q = \pi/4$ have been used for the bounds defined by Equation 110.

with a relative error bound of $8.33 \times 10^{-5}$ for the lower bounded function and $1.44 \times 10^{-4}$ for the upper bounded function. Such accuracy is better than the bounds underpinning the results shown in Figure 17. The sixteenth order approximation, $f_{4,16}$, based on four equal sub-intervals, specified in Appendix 4 and when used with a transition point $x_o = 7.1544$, leads to the functional bounds

$$\frac{f_{4,16}(x)}{1+\varepsilon_B} \leq \text{erf}(x) \leq \frac{f_{4,16}(x)}{1-\varepsilon_B}, \qquad \varepsilon_B = 4.82 \times 10^{-16}, x > 0, \tag{115}$$

with a relative error bound of less than $9.64 \times 10^{-16}$ for the lower bounded function and $9.32 \times 10^{-16}$ for the upper bounded function.

### 7.4    New Series for Error Function

Consider the exact results

$$\text{erf}(x) = f_n(x) + \varepsilon_n(x), \qquad n \in \{0, 1, 2, \ldots\}, \tag{116}$$

where $f_n$ is specified by Equation 24 and $\varepsilon_n^{(1)}(x)$ is specified by Equation 25. By utilizing a Taylor series approximation for $\exp(-x^2)$ in $\varepsilon_n^{(1)}(x)$ and then by integrating, an approximation for $\varepsilon_n(x)$ can be established. This leads to new series for the error function.

### Theorem 7.2    New Series for Error Function

Based on zero, first and second order approximations the following series for the error function are valid:

$$\text{erf}(x) = \frac{x}{\sqrt{\pi}} + \frac{x}{\sqrt{\pi}} \cdot e^{-x^2} + \frac{x^3}{\sqrt{\pi}}\left[\frac{1}{3} - \frac{3x^2}{2 \cdot 5} + \frac{5x^4}{6 \cdot 7} - \frac{7x^6}{24 \cdot 9} + \frac{9x^8}{120 \cdot 11} - \ldots + \frac{(-1)^k(2k+1)x^{2k}}{(2k+3)(k+1)!} + \ldots \right] \tag{117}$$

$$\text{erf}(x) = \frac{x}{\sqrt{\pi}} + \frac{x}{\sqrt{\pi}}\left[1 + \frac{x^2}{3}\right]e^{-x^2} + \frac{x^5}{6\sqrt{\pi}} \cdot \left[\frac{1}{5} - \frac{3x^2}{(3/2) \cdot 7} + \frac{5x^4}{4 \cdot 9} - \frac{7x^6}{15 \cdot 11} + \frac{9x^8}{72 \cdot 13} - \ldots + \frac{(-1)^k(2k+1)(k+1)!x^{2k}}{(2k+5)\prod_{r=1}^{k}\left[\sum_{i=1}^{r} 2i+1\right]} + \ldots \right] \tag{118}$$

$$\operatorname{erf}(x) = \frac{x}{\sqrt{\pi}}\left[1 - \frac{x^2}{30}\right] + \frac{x}{\sqrt{\pi}}\left[1 + \frac{11x^2}{30} + \frac{x^4}{15}\right]e^{-x^2} +$$

$$\frac{1}{\sqrt{\pi}} \cdot \frac{x^7}{60} \cdot \left[\frac{1 \cdot 3}{1 \cdot 3 \cdot 7} - \frac{3 \cdot 5 x^2}{2 \cdot 3 \cdot 9} + \frac{5 \cdot 7 x^4}{4 \cdot 5 \cdot 11} - \frac{7 \cdot 9 x^6}{9 \cdot 10 \cdot 13} + \ldots + \frac{(-1)^k (2k+1)(2k+3)(k+1)! x^{2k}}{3(2k+7) \prod_{r=2}^{k+1}\left[2\sum_{i=2}^{r} i\right]} + \ldots\right] \quad (119)$$

Further series, based on higher order approximations, can also be established.

**Proof**

The proof is detailed in Appendix 7.

### 7.4.1  Results

The relative errors associated with the zero and second order series are shown in Figure 18 and Figure 19. Clearly, the relative error improves as the number of terms used in the series expansion increases. The significant improvement in the relative error, for $x$ small, is evident. A comparison with the relative errors associated with Taylor series approximations, as shown in Figure 2, shows the improved performance.

The second order approximation arising from Equation 118, i.e.

$$\operatorname{erf}(x) = \frac{x}{\sqrt{\pi}} + \frac{x}{\sqrt{\pi}}\left[1 + \frac{x^2}{3}\right]e^{-x^2} + \frac{x^5}{6\sqrt{\pi}} \cdot \left[\frac{1}{5} - \frac{2x^2}{7}\right] \quad (120)$$

yields a relative error bound of less than 0.001 for the interval $[0, 0.87]$ and less than 0.01 for the interval $[0, 1.1]$.

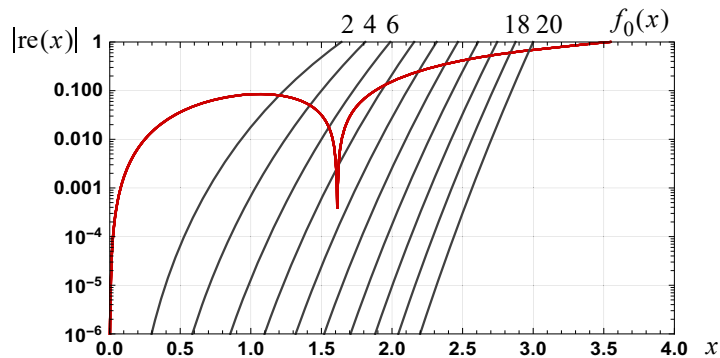

**Figure 18.** Relative error in the approximations $f_0(x)$ and $f_0(x) + \varepsilon_0(x)$ to erf($x$) where the residual function $\varepsilon_0(x)$ is approximated by the stated order.

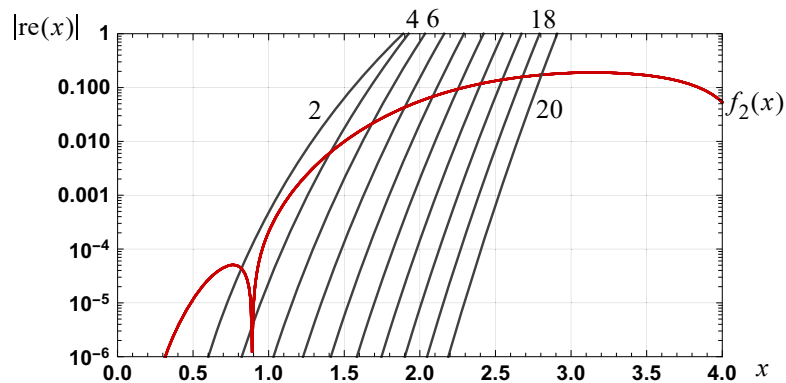

**Figure 19.** Relative error in the approximations $f_2(x)$ and $f_2(x) + \varepsilon_2(x)$ to erf($x$) where the residual function $\varepsilon_2(x)$ is approximated by the stated order.

## 7.5  Complementary Demarcation Functions

Consider a complementary function $e_C$ which is such that

$$e_C^2(x) + \mathrm{erf}^2(x) = 1, \qquad x \geq 0. \tag{121}$$

With the approximation detailed in Theorem 6.1 (and by noting that $\lim_{n \to \infty} p_{n,0} = 1$ - see Equation 96) it is the case that

$$e_C^2(x) = \lim_{n \to \infty}\left[ p_{n,1}(x)e^{-x^2} + p_{n,2}(x)e^{-2x^2} \right] \tag{122}$$

and, thus, $e_C(x)$ can be defined independently of the error function. This function is shown in Figure 20 along with $\mathrm{erf}(x)^2$. These two functions act as complementary demarcation functions for the interval $[0, \infty)$. The transition point is $x_o = 0.74373198514677$ as

$$\mathrm{erf}(x)|_{x = 0.74373198514677} = \frac{1}{\sqrt{2}}. \tag{123}$$

## 7.6  Power and Harmonic Distortion: Erf Modelled Non-linearity

The error function is often used to model nonlinearities and the harmonic distortion created by such a nonlinearity is of interest. Examples include the harmonic distortion in magnetic recording, e.g. Abuelma'atti 1988, Fujiwara, 1980, and the harmonic distortion arising, in a communication context, by a power amplifier, e.g. Taggart 2005. For these cases the interest was in obtaining, with a sinusoidal input signal defined by $a\sin(2\pi f_o t)$, the harmonic distortion created by an error function nonlinearity over the input amplitude range of $[-2, 2]$.

Consider the output signal of a nonlinearity modelled by the error function:

$$y(t) = \mathrm{erf}[a\sin(2\pi f_o t)]. \tag{124}$$

For such a case, the output power is defined according to

$$P = \frac{1}{T} \cdot \int_0^T y^2(t)dt = \frac{1}{T} \cdot \int_0^T \mathrm{erf}^2[a\sin(2\pi f_o t)]dt, \qquad T = 1/f_o, \tag{125}$$

and the output amplitude associated with the $k$th harmonic is

$$\frac{\sqrt{2}}{\sqrt{T}}\int_0^T \mathrm{erf}[a\sin(2\pi f_o t)]\sin(2\pi k f_o t)dt. \tag{126}$$

To determine an analytical approximation to the output power, the approximations stated in Theorem 6.1 lead to relatively simple expressions. Consider the third order approximation, as specified by Equation 93, which has a relative error bound of $2.03 \times 10^{-4}$ for the positive real line. For such a case, the output signal is approximated according to

$$y_3(t) = \frac{1}{\sqrt{\pi}} \left[ \begin{array}{l} \dfrac{22}{7} - \dfrac{40 e^{-[a\sin(2\pi f_o t)]^2}}{21}\left[1 - \dfrac{[a\sin(2\pi f_o t)]^2}{20}\right] - \\ \dfrac{26 e^{-2[a\sin(2\pi f_o t)]^2}}{21}\left[1 + \dfrac{10[a\sin(2\pi f_o t)]^2}{26} + \dfrac{[a\sin(2\pi f_o t)]^4}{13} + \dfrac{[a\sin(2\pi f_o t)]^6}{130}\right] \end{array} \right]^{1/2} \quad (127)$$

and is shown in Figure 21.

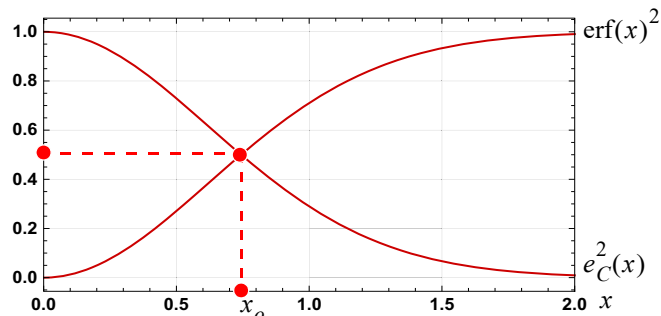

**Figure 20.** Graph of the signals $e_C^2(x)$ and $\mathrm{erf}(x)^2$.

The power in $y_3$ can be readily be determined (e.g. via use of Mathematica) and it then follows that an approximation to the true power is

$$P(a) \approx \frac{22}{7\pi} - \frac{40}{21\pi} I_0\left[\frac{a^2}{2}\right]\left[1 - \frac{a^2}{40}\right] e^{-a^2/2} - \frac{26}{21\pi} I_0[a^2]\left[1 + \frac{5a^2}{26} + \frac{41 a^4}{1040} + \frac{a^6}{260}\right] e^{-a^2} - \frac{a^2}{21\pi} I_1\left[\frac{a^2}{2}\right] e^{-a^2/2} + \frac{37 a^2}{140\pi} I_1[a^2]\left[1 + \frac{43 a^2}{222} + \frac{2 a^4}{111}\right] e^{-a^2} \quad (128)$$

where $I_0$ and $I_1$, respectively, are the zero and first order Bessel functions of the first kind. The variation in output power is shown in Figure 22.

### 7.6.1   Harmonic Distortion

To establish analytical approximations for the harmonic distortion, the functional forms detailed in Theorem 6.1 are not suitable. However, the functional forms detailed in Theorem 2.1 do lead to analytical approximations which are valid over a restricted domain. Consider, a fourth order spline approximation, as specified by Equation 36, which approximates the error function over the range $[-2, 2]$ with a relative error bound that is better than $0.001$ and leads to the approximation

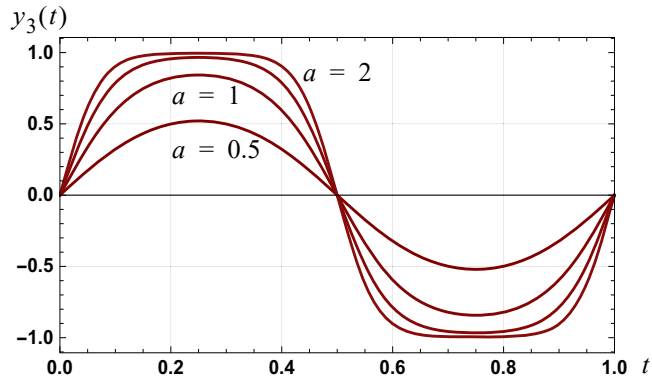

**Figure 21.** Graph of $y_3(t)$ for the case of $f_o = 1$ and for amplitudes of $a = 0.5$, $a = 1$, $a = 1.5$ and $a = 2$.

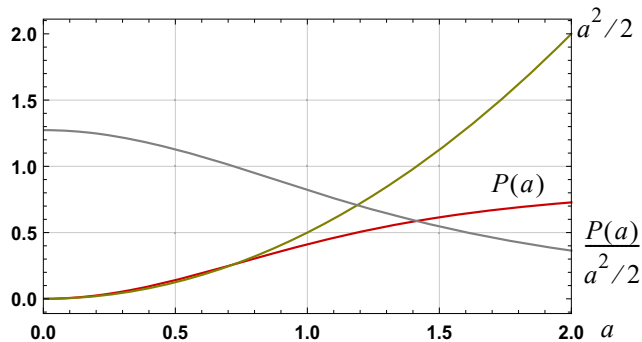

**Figure 22.** Graph of the input power, output power and ratio of output power to input power as the amplitude of the input signal varies.

$$y_4(t) = \frac{a\sin(2\pi f_o t)}{\sqrt{\pi}}\left[1 - \frac{a^2\sin(2\pi f_o t)^2}{18} + \frac{a^4\sin(2\pi f_o t)^4}{1260}\right] +$$

$$\frac{a\sin(2\pi f_o t)}{\sqrt{\pi}}\left[1 + \frac{7a^2\sin(2\pi f_o t)^2}{18} + \frac{37a^4\sin(2\pi f_o t)^4}{420} + \frac{4a^6\sin(2\pi f_o t)^6}{315} + \frac{a^8\sin(2\pi f_o t)^8}{945}\right]e^{-a^2\sin(2\pi f_o t)^2} \quad (129)$$

The amplitude of the $k$th harmonic in such a signal is given by

$$c_{4,k} = \frac{\sqrt{2}}{\sqrt{T}}\int_0^T y_4[a\sin(2\pi f_o t)]\sin(2\pi k f_o t)dt = \sqrt{2T}\int_0^1 y_4[a\sin(2\pi\lambda)]\sin(2\pi k\lambda)d\lambda \quad (130)$$

where the change of variable $\lambda = f_o t$ has been used. The first, third, fifth and seventh harmonic levels are:

$$\frac{c_{4,1}}{\sqrt{T}} = \frac{\sqrt{2}a}{2\sqrt{\pi}}\left[1 - \frac{a^2}{24} + \frac{a^4}{2016}\right] + \frac{\sqrt{2}a}{2\sqrt{\pi}} \cdot I_0\left[\frac{a^2}{2}\right]\left[1 + \frac{11a^2}{24} + \frac{11a^4}{105} + \frac{a^6}{70} + \frac{a^8}{945}\right]e^{-a^2/2} -$$

$$\frac{5\sqrt{2}a}{6\sqrt{\pi}} \cdot I_1\left[\frac{a^2}{2}\right]\left[1 + \frac{1481a^2}{4200} + \frac{38a^4}{525} + \frac{29a^6}{3150} + \frac{a^8}{1575}\right]e^{-a^2/2} \quad (131)$$

$$\frac{c_{4,3}}{\sqrt{T}} = \frac{\sqrt{2}a^3}{144\sqrt{\pi}}\left[1 - \frac{a^2}{56}\right] - \frac{115\sqrt{2}a}{84\sqrt{\pi}} \cdot I_0\left[\frac{a^2}{2}\right]\left[1 + \frac{403a^2}{1380} + \frac{6a^4}{115} + \frac{31a^6}{5175} + \frac{2a^8}{5175}\right]e^{-a^2/2} +$$

$$\frac{115\sqrt{2}}{21a\sqrt{\pi}} \cdot I_1\left[\frac{a^2}{2}\right]\left[1 + \frac{8a^2}{23} + \frac{163a^4}{1840} + \frac{76a^6}{5175} + \frac{11a^8}{6900} + \frac{a^{10}}{10350}\right]e^{-a^2/2} \quad (132)$$

$$\frac{c_{4,5}}{\sqrt{T}} = \frac{\sqrt{2}a^5}{40320\sqrt{\pi}} + \frac{262\sqrt{2}}{15a\sqrt{\pi}} \cdot I_0\left[\frac{a^2}{2}\right]\left[1 + \frac{1943a^2}{7336} + \frac{1485a^4}{29344} + \frac{73a^6}{11004} + \frac{13a^8}{22008} + \frac{a^{10}}{33012}\right]e^{-a^2/2} -$$

$$\frac{1048\sqrt{2}}{15a^3\sqrt{\pi}} \cdot I_1\left[\frac{a^2}{2}\right]\left[1 + \frac{1943a^2}{7336} + \frac{1201a^4}{14672} + \frac{5125a^6}{352128} + \frac{5a^8}{2751} + \frac{41a^{10}}{264096} + \frac{a^{12}}{132048}\right]e^{-a^2/2} \quad (133)$$

$$\frac{c_{4,7}}{\sqrt{T}} = \frac{-6784\sqrt{2}}{21a^3\sqrt{\pi}} \cdot I_0\left[\frac{a^2}{2}\right]\left[1 + \frac{779a^2}{3392} + \frac{631a^4}{13568} + \frac{2047a^6}{325632} + \frac{25a^8}{40704} + \frac{17a^{10}}{407040} + \frac{a^{12}}{610560}\right]e^{-a^2/2} -$$

$$\frac{27136\sqrt{2}}{21a^5\sqrt{\pi}} \cdot I_1\left[\frac{a^2}{2}\right]\left[1 + \frac{779a^2}{3392} + \frac{1055a^4}{13568} + \frac{137a^6}{10176} + \frac{2269a^8}{1302528} + \frac{67a^{10}}{407040} + \frac{a^{12}}{92160} + \frac{a^{14}}{2442240}\right]e^{-a^2/2} \quad (134)$$

The variation, with the input signal amplitude, of the harmonic distortion, as defined by $c_{4,k}^2/c_{4,1}^2$, is shown in Figure 23.

### 7.7    Linear Filtering of a Error Function Step Signal

Consider the case of a practical step input signal that is modelled by the error function, $\text{erf}(t/\gamma)$ and the case where such a signal is input to a 2nd order linear filter with a transfer function defined by

$$H(s) = \frac{1}{\left[1 + \frac{s}{2\pi f_p}\right]^2} \quad \Leftrightarrow \quad h(t) = \frac{te^{-t/\tau}}{\tau^2} \cdot u(t), \quad \tau = \frac{1}{2\pi f_p}. \quad (135)$$

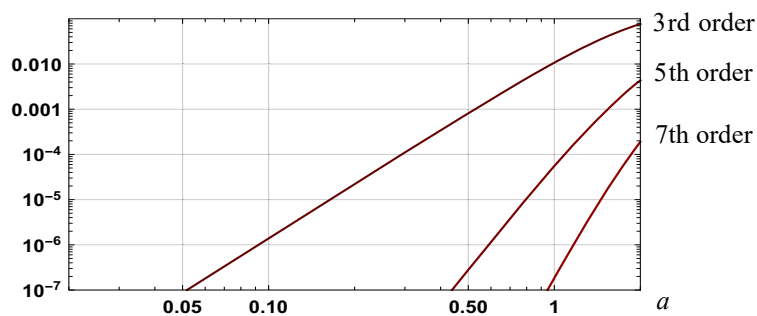

**Figure 23.** Graph of the variation of harmonic distortion with amplitude.

### Theorem 7.3   Linear Filtering of an Error Function Signal

The output of a second order linear filter, defined by Equation 135, to an error function input signal, defined by $\mathrm{erf}(t/\gamma)$, is

$$y(t) = \mathrm{erf}\left[\frac{t}{\gamma}\right]u(t) + \frac{e^{-t/\tau}}{\tau} \cdot \left\{\left[\frac{\gamma^2}{2\tau} - (t+\tau)\right]\exp\left[\frac{\gamma^2}{4\tau^2}\right]\left[\mathrm{erf}\left[\frac{\gamma}{2\tau}\right] - \mathrm{erf}\left[\frac{\gamma}{2\tau} - \frac{t}{\gamma}\right]\right] - \frac{\gamma}{\sqrt{\pi}}\exp\left[\frac{\gamma^2}{4\tau^2}\right]\exp\left[-\left[\frac{t}{\gamma} - \frac{\gamma}{2\tau}\right]^2\right] + \frac{\gamma}{\sqrt{\pi}}\right\}u(t) \quad (136)$$

and can be approximated by the $n$th order signal

$$y_n(t) = f_n\left[\frac{t}{\gamma}\right]u(t) + \frac{e^{-t/\tau}}{\tau} \cdot \left\{\left[\frac{\gamma^2}{2\tau} - (t+\tau)\right]\exp\left[\frac{\gamma^2}{4\tau^2}\right]\left[f_n\left[\frac{\gamma}{2\tau}\right] - f_n\left[\frac{\gamma}{2\tau} - \frac{t}{\gamma}\right]\right] - \frac{\gamma}{\sqrt{\pi}}\exp\left[\frac{\gamma^2}{4\tau^2}\right]\exp\left[-\left[\frac{t}{\gamma} - \frac{\gamma}{2\tau}\right]^2\right] + \frac{\gamma}{\sqrt{\pi}}\right\}u(t) \quad (137)$$

where $f_n$ is defined by one of the approximations detailed in Theorem 3.1, Theorem 4.1, Theorem 5.1 or Theorem 6.1. It is the case that

$$\lim_{n \to \infty} y_n(t) = y(t), \quad t \in (0, \infty). \quad (138)$$

### Proof

The proof is detailed in Appendix 8.

#### 7.7.1    Results

For an input signal $\mathrm{erf}(t/\gamma)u(t)$, $\gamma = 1/2$, input into a second order linear filter with $f_p = 1$, the output signal is shown in Figure 24. The relative errors in the approximations to the output signal are shown in Figure 25 for the case of approximations as specified by Equation 56 and with the use of optimum transition points.

### 7.8    Extension to Complex Case

By definition, the error function, for the general complex case, is defined according to

$$\mathrm{erf}(z) = \frac{2}{\sqrt{\pi}}\int_\gamma e^{-\lambda^2}d\lambda, \quad z \in \mathbf{C}, \quad (139)$$

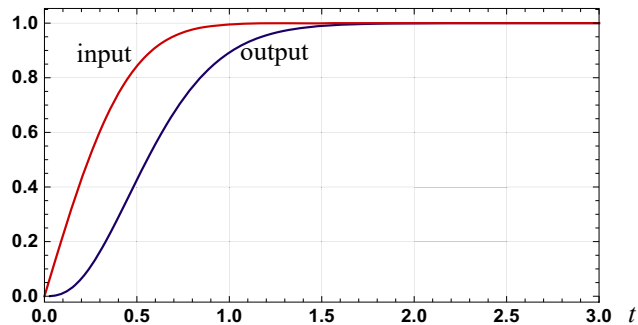

**Figure 24.** Graph of the input signal $\mathrm{erf}(t/\gamma)$, $\gamma = 1/2$ and the corresponding approximation to the output of a second order linear filter with $f_p = 1$, $\tau = 1/2\pi$.

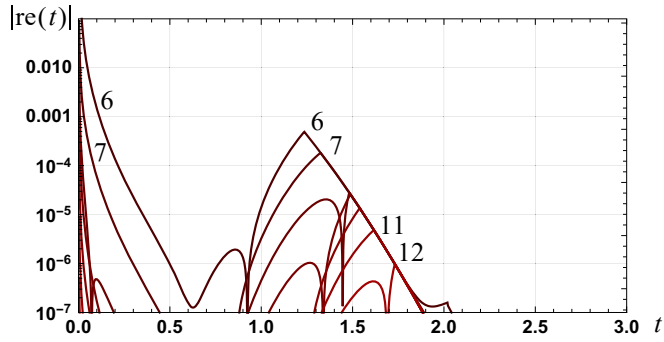

**Figure 25.** Graph of the relative errors associated with the output signal, shown in Figure 24, for approximations to the error function (Equation 56) of orders six to twelve which utilize optimum transition points.

where the path $\gamma$ is between the points zero and $z$ and is arbitrary. For the case of $z = x + jy$, and a path along the $x$ axis to the point $(x, 0)$ and then to the point $z$, the error function is defined according to (Salzer 1951, eqn. 5)

$$\begin{aligned}\text{erf}(x+jy) &= \frac{2}{\sqrt{\pi}}\int_0^x e^{-\lambda^2}d\lambda + \frac{2j}{\sqrt{\pi}}\int_0^y e^{-(x+j\lambda)^2}d\lambda \\ &= \text{erf}(x) + \frac{2e^{-x^2}}{\sqrt{\pi}}\cdot\int_0^y e^{\lambda^2}\sin(2x\lambda)d\lambda + \frac{2je^{-x^2}}{\sqrt{\pi}}\cdot\int_0^y e^{\lambda^2}\cos(2x\lambda)d\lambda\end{aligned} \tag{140}$$

Explicit approximations for $\text{erf}(x+jy)$ then arise when integrable approximations for the two dimensional surfaces $\exp(y^2)\sin(2xy)$ and $\exp(y^2)\cos(2xy)$ over $[0, x] \times [0, y]$ are available. Naturally, significant existing research exists, e.g. Salzer 1951, Abrarov, 2013.

### 7.9     Approximation for the Inverse Error Function

There are many applications where the inverse error function is required and accurate approximations for this function is of interest. From the research underpinning this paper, the authors view is that finding approximations to the inverse error function is best treated directly and as a separate problem, rather than approaching it via finding the inverse of an approximation to the error function.

## 8     Conclusion

This paper has detailed analytical approximations for the real case of the error function, underpinned by a spline based integral approximation, which have significantly better convergence than the default Taylor series. The original approximations can be improved by utilizing the approximation $\text{erf}(x) \approx 1$ for $x > x_o$, with $x_o$ being dependent on the order of approximation. The fourth order approximations arising from Theorem 2.1 and Theorem 2.3, with respective transition points of $x_o = 2.3715$ and $x_o = 2.6305$, achieve relative error bounds over the interval $[0, \infty]$, respectively, of $1.03 \times 10^{-3}$ and $2.28 \times 10^{-4}$. The respective sixteenth order approximations, with $x_o = 3.9025$ and $x_o = 4.101$, have relative error bounds of $3.44 \times 10^{-8}$ and $6.66 \times 10^{-9}$.

Further improvements were detailed via two generalizations. The first was based on utilizing integral approximations for each of the $m$ equally spaced sub-intervals in the required interval of integration. The second, was based on utilizing a fixed sub-interval within the interval of integration, with a known tabulated area, and then utilizing an integral approximation over the remainder of the interval. Both generalizations lead to significantly improved accuracy. For example, a fourth order approximation based on four sub-intervals, with $x_o = 3.7208$, achieves a relative error bound of $1.43 \times 10^{-7}$ over the interval $[0, \infty]$. A sixteenth order approximation, with $x_o = 6.3726$, has a relative error bound of $2.01 \times 10^{-19}$.

Finally, it was shown that a custom feedback system, with inputs defined by either the original error function approximations or approximations based on the use of sub-intervals, leads to analytical approximations with improved accuracy and which are valid over the positive real line without utilizing the approximation $\text{erf}(x) \approx 1$ for $x$ suitably large. The original fourth order error function approximation yields an approximation with a rela-tive error bound of $1.82 \times 10^{-5}$ over the interval $[0, \infty]$. The original sixteenth order approximation yields an approximation with a relative error bound of $1.68 \times 10^{-14}$.

Applications of the approximations were detailed and these include, first, approximations to achieve the specified error bounds of $10^{-4}$, $10^{-6}$, $10^{-10}$ and $10^{-16}$ over the positive real line. Second, the definitions of functions that are upper and lower bounds, of arbitrary accuracy, for the error function. Third, new series for the error function. Fourth, new sequences of approximations for $\exp(-x^2)$ which have significantly higher convergence properties that a Taylor series approximation. Fifth, a complementary demarcation function satisfying the constraint $e_C^2(x) + \mathrm{erf}^2(x) = 1$ was defined. Sixth, arbitrarily accurate approximations for the power and harmonic distortion for a sinusoidal signal subject to an error function nonlinearity. Seventh, approximate expressions for the linear filtering of a step signal that is modelled by the error function.

**Acknowledgement**   The support of Prof. A. Zoubir, SPG, Technische Universität Darmstadt, Darmstadt, Germany, who hosted a visit where the research for, and writing of, this paper was completed, is gratefully acknowledged. An anonymous reviewer provided a comment, which after consideration, led to the improved approximations detailed in Theorem 2.3.

# Appendix 1:    Proof of Theorem 2.1

Consider $f(x) = \exp(-x^2)$. Successive differentiation of this function leads to the iterative formula

$$f^{(k)}(x) = p(k, x)\exp(-x^2) \tag{141}$$

where

$$p(k, x) = p^{(1)}(k-1, x) - 2xp(k-1, x), \qquad p(0, x) = 1. \tag{142}$$

It then follows from Equation 20 that

$$\frac{2}{\sqrt{\pi}} \cdot \int_\alpha^x \exp(-\lambda^2)d\lambda \approx \frac{2}{\sqrt{\pi}} \cdot \sum_{k=0}^{n} c_{n,k}(x-\alpha)^{k+1}[f^{(k)}(\alpha) + (-1)^k f^{(k)}(x)]$$

$$= \frac{2}{\sqrt{\pi}} \cdot \sum_{k=0}^{n} c_{n,k}(x-\alpha)^{k+1}[p(k, \alpha)\exp(-\alpha^2) + (-1)^k p(k, x)\exp(-x^2)] \tag{143}$$

The result for the case of $\alpha = 0$ then yields the $n$th order approximation for the error function:

$$f_n(x) = \frac{2}{\sqrt{\pi}} \cdot \sum_{k=0}^{n} c_{n,k} x^{k+1}\left[p(k, 0) + (-1)^k p(k, x)e^{-x^2}\right] \tag{144}$$

To determine $\varepsilon_n^{(1)}(x)$ consider the equality $\mathrm{erf}(x) = f_n(x) + \varepsilon_n(x)$. Differentiation yields

$$\varepsilon_n^{(1)}(x) = \frac{2e^{-x^2}}{\sqrt{\pi}} - \frac{2}{\sqrt{\pi}} \cdot \sum_{k=0}^{n} c_{n,k}(k+1)x^k\left[p(k, 0) + (-1)^k p(k, x)e^{-x^2}\right] -$$

$$\frac{2e^{-x^2}}{\sqrt{\pi}} \cdot \sum_{k=0}^{n} c_{n,k} x^{k+1}(-1)^k[p^{(1)}(k, x) - 2xp(k, x)] \tag{145}$$

and the required result:

$$\varepsilon_n^{(1)}(x) = \frac{2e^{-x^2}}{\sqrt{\pi}} - \frac{2}{\sqrt{\pi}} \cdot \sum_{k=0}^{n} c_{n,k}(k+1)x^k p(k, 0) -$$

$$\frac{2e^{-x^2}}{\sqrt{\pi}} \cdot \sum_{k=0}^{n} c_{n,k} x^k(-1)^k[(k+1-2x^2)p(k, x) + xp^{(1)}(k, x)] \tag{146}$$

## Appendix 2:    Mathematica Code (Version 12.3) for Figure 2 Results

ClearAll["Global`*"];

**(* Define Parameters *)**

orderMin=0;     orderMax=10;

xMin=0;   xMax=4;   yMax=1;   yMin=10$^{-6}$;
nPoint=1000;   resolution=(xMax-xMin)/nPoint;

aRatio=6/11;    iSize={310,170};
font1={"Courier",FontSize->9,FontWeight->"Bold"};
pStyleS={RGBColor[1/4,0,0],PointSize[0.005],Thickness[0.004]};
pStyleT={GrayLevel[0.],PointSize[0.005],Thickness[0.004]};

**(* Define Integrand *)**

f[x_]:=$\frac{2}{\sqrt{\pi}}$ *Exp[-x$^2$];     df[x_,i_]:=D[f[tx],{tx,i}]/.tx->x;

Plot[f[x],{x,xMin,xMax},PlotRange->{{xMin,xMax},All},Frame->True,GridLines->Automatic,FrameTicks->Automatic,AspectRatio->aRatio,ImageSize->iSize,PlotStyle->pStyleS,BaseStyle->font1,PlotLabel->"Integrand"]

**(* Define Spline Approx: Two Forms *)**

p[x_,0]=1;
p[x_,k_]:=p[x,k]=Expand[(D[p[tx,k-1],{tx,1}]/.tx->x) - 2*x*p[x,k-1]]

Print[" "]
Do[
 Print["order = ",i," Definition for p[x,i] = ",p[x,i] ],
  {i,0,6,1}]
Print[" "]

$$\text{erfApprox}[n\_,x1\_,x2\_]:= \sum_{k=0}^{n} \frac{n!}{(n-k)!\cdot(k+1)!} \cdot \frac{(2\cdot n+1-k)!}{2\cdot(2\cdot n+1)!} \cdot (x2-x1)^{k+1} \cdot (df[x1,k] + (-1)^k \cdot df[x2,k])$$

erfApproxAlt[n_,x1_,x2_]:=

$$\frac{2}{\sqrt{\pi}} \cdot \sum_{k=0}^{n} \frac{n!}{(n-k)!\cdot(k+1)!} \cdot \frac{(2\cdot n+1-k)!}{2\cdot(2\cdot n+1)!} \cdot (x2-x1)^{k+1} \cdot \left(p(x1,k)\cdot e^{-x1^2} + (-1)^k \cdot p(x2,k)\cdot e^{-x2^2}\right)$$

re[order_,x_]:=1-erfApprox[order,0,x]/Erf[x]

**(* Explicit Spline Approx. Check on Approx. *)**

$$fE[0,x\_]:\ =\ \frac{x}{\sqrt{\pi}} + \frac{x}{\sqrt{\pi}} \cdot e^{-x^2}$$

$$fE[1,x\_]:\ =\ \frac{x}{\sqrt{\pi}} + \frac{x}{\sqrt{\pi}} \cdot \left(1 + \frac{x^2}{3}\right) \cdot e^{-x^2}$$

$$fE[2,x\_]:\ =\ \frac{x}{\sqrt{\pi}} \cdot \left(1 - \frac{x^2}{30}\right) + \frac{x}{\sqrt{\pi}} \cdot \left(1 + \frac{11\cdot x^2}{30} + \frac{x^4}{15}\right) \cdot e^{-x^2}$$

$$fE[3,x\_]:\ =\ \frac{x}{\sqrt{\pi}} \cdot \left(1 - \frac{x^2}{21}\right) + \frac{x}{\sqrt{\pi}} \cdot \left(1 + \frac{8\cdot x^2}{21} + \frac{17\cdot x^4}{210} + \frac{x^6}{105}\right) \cdot e^{-x^2}$$

$$fE[4,x\_]:\ =\ \frac{x}{\sqrt{\pi}} \cdot \left(1 - \frac{x^2}{18} + \frac{x^4}{1260}\right) + \frac{x}{\sqrt{\pi}} \cdot \left(1 + \frac{7\cdot x^2}{18} + \frac{37\cdot x^4}{420} + \frac{4\cdot x^6}{315} + \frac{x^8}{945}\right) \cdot e^{-x^2}$$

$$fE[5,x\_]:\ =\ \frac{x}{\sqrt{\pi}} \cdot \left(1 - \frac{2\cdot x^2}{33} + \frac{x^4}{660}\right) + \frac{x}{\sqrt{\pi}} \cdot \left(1 + \frac{13\cdot x^2}{33} + \frac{61\cdot x^4}{660} + \frac{67\cdot x^6}{4620} + \frac{16\cdot x^8}{10395} + \frac{x^{10}}{10395}\right) \cdot e^{-x^2}$$

Print[" "]

```
Do[
        Print["order  =  ",i,"  Check  on  Spline  Approx:  ",Simplify[erfApprox[i,0,x]-fE[i,x]],"
            ",Simplify[erfApproxAlt[i,0,x] - fE[i,x]]  ],
{i,0,5,1}]

Print[" "]
```

**(* RE for Spline Approx. *)**

```
Do[
        data=Table[Abs[N[re[order,xMin+i*resolution]]],{i,1,nPoint,1}];
        dataG=Table[{xMin+i*resolution,Abs[N[re[order,xMin+i*resolution]]]}, {i,1,nPoint,1}];

        maxRE[order]=Max[data];

        Print[pRE[order]=ListLogPlot[dataG, PlotRange->{{xMin,xMax},{yMin,yMax}}, Joined->True, Frame-
            >True, Grid Lines->Automatic, FrameTicks->Automatic, AspectRatio->aRatio, ImageSize->iSize,
            PlotStyle->pStyleS, BaseStyle->font1, PlotLabel->{"RE: Erf Approx: Order = ",order," Max RE =
            ",NumberForm[maxRE[order],5]}]   ],

{order,orderMin,orderMax,1}];

tableS=Table[pRE[order],{order,orderMin,orderMax,1}];
eS=Show[tableS,PlotLabel->" "]
```

**(* Table of RE bounds for [0,4] *)**

```
dataRE=Table[{order,NumberForm[maxRE[order],5]},{order,orderMin,orderMax,1}];
dataRE=Prepend[dataRE,{SplineOrder,MaxRE}];
Print[StyleForm[TableForm[dataRE],FontSize->10]];
```

**(* Taylor Series *)**

```
orderMinT=1;       orderMaxT=15;

resolutionT=(xMax-xMin)/100;
```

$$\text{taylorSeries}[\text{order}\_,x\_]:= \frac{2}{\sqrt{\pi}} \cdot \sum_{k=1}^{\text{order}} \frac{1+(-1)^{k+1}}{2} \cdot \frac{(-1)^{(k-1)/2} \cdot x^k}{k \cdot \text{Factorial}\left[\frac{k-1}{2}\right]}$$

```
reT[order_,x_]:=1-taylorSeries[order,x]/Erf[x]

Print[" "]
Do[
   Print["Taylor series: order = ",i,"   ",Expand[taylorSeries[i,x]] ],
{i,1,9,2}];

Print[" "];

Do[
        tableRET=Table[{i*resolutionT,Abs[N[reT[orderT,i*resolutionT]]]},{i,1,nPoint,1}];

        pRET[orderT]=ListLogPlot[tableRET, PlotRange->{{xMin,xMax},{yMin,yMax}}, Joined->False,
        Frame->True, GridLines->Automatic, FrameTicks->Automatic, AspectRatio->aRatio, ImageSize->iSize,
        PlotStyle->pStyleT, BaseStyle->font1, PlotLabel->{"RE in Erf Approx. Order = ",order}],

{orderT,orderMinT,orderMaxT,2}]

tableT=Table[pRET[orderT],{orderT,orderMinT,orderMaxT,2}];
eT=Show[tableT,PlotLabel->"Taylor Series Approx.  "]
```

**(* Combined Graph: Taylor Series + Spline Approx *)**

```
e1=Show[eS,eT,PlotLabel->" "]
```

## Appendix 3:   Proof of Theorem 2.2

The use of a Taylor series expansion for $\exp(-x^2)$ in the definitions of $\varepsilon_0^{(1)}(x)$ and $\varepsilon_1^{(1)}(x)$, as defined by Equation 40 and Equation 41, yields:

$$\varepsilon_0^{(1)}(x) = \frac{x^2}{\sqrt{\pi}} \cdot \left[1 - \frac{3x^2}{2} + \frac{5x^4}{6} - \frac{7x^6}{24} + \frac{3x^8}{40} - \frac{11x^{10}}{720} + \frac{13x^{12}}{5040} - \frac{x^{14}}{2688} + \ldots \right]$$

$$= \frac{x^2}{\sqrt{\pi}} \cdot [c_{0,0} - c_{0,1}x^2 + \ldots + (-1)^k c_{0,k} x^{2k} + \ldots], \qquad c_{0,k} = \frac{2}{k!} - \frac{1}{(k+1)!}, \quad k \geq 0 \tag{147}$$

$$\varepsilon_1^{(1)}(x) = \frac{x^4}{6\sqrt{\pi}} \cdot \left[1 - 2x^2 + \frac{5x^4}{4} - \frac{7x^6}{16} + \frac{x^8}{8} - \frac{11x^{10}}{420} \ldots \right] \tag{148}$$

From a consideration of $\varepsilon_0^{(1)}(x)$ and $\varepsilon_1^{(1)}(x)$, as well as higher order residual functions, it can be readily seen that the polynomial terms of order zero to $2n+1$ in $\varepsilon_n^{(1)}(x)$ have coefficients of zero. It then follows that $\varepsilon_n^{(1)}(x)$ can be written as

$$\varepsilon_n^{(1)}(x) = \frac{1}{\sqrt{\pi}} \cdot \frac{x^{2n+2}}{x_{n,0}} \cdot g_n(x), \qquad n \in \{0, 1, 2, \ldots\}, \tag{149}$$

where

$$g_n(x) = 1 - d_{n,1}x^2 + d_{n,2}x^4 - \ldots + (-1)^k d_{n,k} x^{2k} + \ldots \tag{150}$$

and where it can readily be shown that

$$x_{n,0} = 2^n \prod_{i=0}^{n} (2i+1) \geq 2^n n!, \qquad n \in \{0, 1, 2, \ldots\}. \tag{151}$$

Graphs of the magnitude of the residual functions, $\varepsilon_n^{(1)}(x)$, for orders zero, two, four, six and eight, are shown in Figure 26. The magnitude of the functions defined by $g_n$, for the same orders, are shown in Figure 27 and it is evident that

$$|g_n(x)| \leq k_o, \qquad x \geq 0, n \in \{0, 1, 2, \ldots\}, \tag{152}$$

for a fixed constant $k_o$ which is of the order of unity. Hence:

$$|\varepsilon_n^{(1)}(x)| \leq \frac{k_o}{\sqrt{\pi}} \cdot \frac{x^{2n+2}}{x_{n,0}}, \qquad x \geq 0, n \in \{0, 1, 2, \ldots\}. \tag{153}$$

Further, as $x_{n,0} \geq 2^n n!$, it follows, for all fixed values of $x$, that

$$\lim_{n \to \infty} \varepsilon_n^{(1)}(x) = 0, \qquad x \geq 0. \tag{154}$$

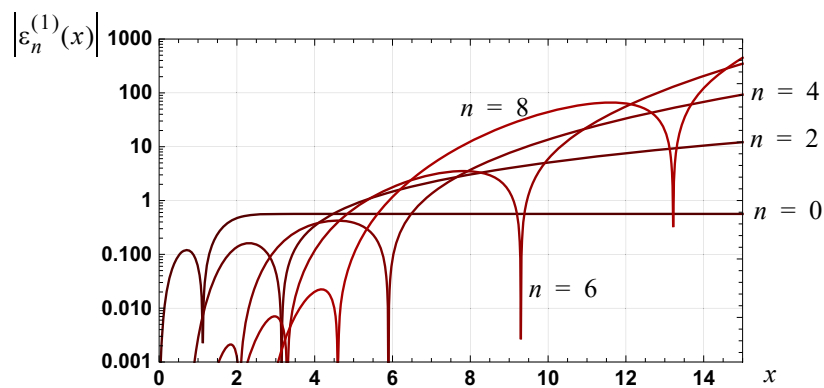

**Figure 26.** Graphs of $|\varepsilon_n^{(1)}(x)|$ for orders zero, two, four, six and eight.

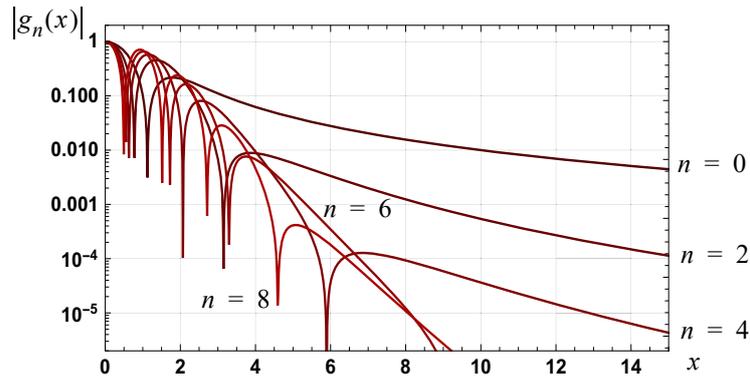

**Figure 27.** Graphs of $|g_n(x)|$ for orders zero, two, four, six and eight.

The convergence is not uniform.

It then follows, for all fixed values of $x$, that there exists an order of approximation, $n$, such that the error in the approximation $\varepsilon_n^{(1)}(x)$ can be made arbitrarily small, i.e. $\forall \varepsilon_o > 0$ there exists a number $N(x)$ such that

$$\left|\varepsilon_n^{(1)}(x)\right| < \varepsilon_o, \qquad \forall n > N(x). \tag{155}$$

In general, $N(x)$ increases with $x$. Thus, $\forall \varepsilon_o > 0$ there exists a number $N_{x_o}(x_o)$ such that

$$\left|\varepsilon_n^{(1)}(x)\right| < \varepsilon_o, \qquad x \in [0, x_o], \forall n > N_{x_o}(x_o). \tag{156}$$

Finally, as $\varepsilon_n^{(1)}(0) = 0$, for all $n$, it then follows, for $x$ fixed, that

$$|\varepsilon_n(x)| = \left|\int_0^x \varepsilon_n^{(1)}(\lambda) d\lambda\right| < \varepsilon_o x, \qquad \forall n > N_x(x), \tag{157}$$

which proves convergence.

## Appendix 4: Fourth Order Spline Approximation - Sixteen Sub-interval Case

Consistent with Theorem 4.1, a fourth order spline approximation, which utilizes sixteen sub-intervals, is

$$f_{4, 16}(x) = \frac{x}{16\sqrt{\pi}}\left[1 - \frac{16x^2}{73728} + \frac{16x^4}{1,321,205,760}\right] +$$

$$\frac{x}{8\sqrt{\pi}} \cdot \exp\left[\frac{-x^2}{256}\right]\left[1 - \frac{x^2}{4608} + \frac{47x^4}{27,525,120} - \frac{x^6}{5,284,823,040} + \frac{x^8}{4,058,744,094,720}\right] +$$

$$\frac{x}{8\sqrt{\pi}} \cdot \exp\left[\frac{-x^2}{64}\right]\left[1 - \frac{x^2}{4608} + \frac{187x^4}{27,525,120} - \frac{x^6}{1,321,205,760} + \frac{x^8}{253,671,505,920}\right] +$$

$$\frac{x}{8\sqrt{\pi}} \cdot \exp\left[\frac{-9x^2}{256}\right]\left[1 - \frac{x^2}{4608} + \frac{1261x^4}{82,575,360} - \frac{x^6}{587,202,560} + \frac{3x^8}{150,323,855,360}\right] +$$

$$\frac{x}{8\sqrt{\pi}} \cdot \exp\left[\frac{-x^2}{16}\right]\left[1 - \frac{x^2}{4608} + \frac{249x^4}{9,175,040} - \frac{x^6}{330,301,440} + \frac{x^8}{15,854,469,120}\right] +$$

$$\frac{x}{8\sqrt{\pi}} \cdot \exp\left[\frac{-25x^2}{256}\right]\left[1 - \frac{x^2}{4608} + \frac{389x^4}{9,175,040} - \frac{5x^6}{1,056,964,608} + \frac{125x^8}{811,748,818,944}\right] +$$

$$\frac{x}{8\sqrt{\pi}} \cdot \exp\left[\frac{-9x^2}{64}\right]\left[1 - \frac{x^2}{4608} + \frac{5041x^4}{82,575,360} - \frac{x^6}{146,800,640} + \frac{3x^8}{9,395,240,960}\right] +$$

$$\frac{x}{8\sqrt{\pi}} \cdot \exp\left[\frac{-49x^2}{256}\right]\left[1 - \frac{x^2}{4608} + \frac{2287x^4}{27,525,120} - \frac{7x^6}{754,974,720} + \frac{343x^8}{579,820,584,960}\right] +$$

$$\frac{x}{8\sqrt{\pi}} \cdot \exp\left[\frac{-x^2}{4}\right]\left[1 - \frac{x^2}{4608} + \frac{2987x^4}{27,525,120} - \frac{x^6}{82,575,360} + \frac{x^8}{990,904,320}\right] + \quad (158)$$

$$\frac{x}{8\sqrt{\pi}} \cdot \exp\left[\frac{-81x^2}{256}\right] \cdot \left[1 - \frac{x^2}{4608} + \frac{11341x^4}{82,575,360} - \frac{9x^6}{587,202,560} + \frac{243x^8}{150,323,855,360}\right] +$$

$$\frac{x}{8\sqrt{\pi}} \cdot \exp\left[\frac{-25x^2}{64}\right]\left[1 - \frac{x^2}{4608} + \frac{4667x^4}{27,525,120} - \frac{5x^6}{264,241,152} + \frac{125x^8}{50,734,301,184}\right] +$$

$$\frac{x}{8\sqrt{\pi}} \cdot \exp\left[\frac{-121x^2}{256}\right]\left[1 - \frac{x^2}{4608} + \frac{5647x^4}{27,525,120} - \frac{121x^6}{5,284,823,040} + \frac{14641x^8}{4,058,744,094,720}\right] +$$

$$\frac{x}{8\sqrt{\pi}} \cdot \exp\left[\frac{-9x^2}{16}\right]\left[1 - \frac{x^2}{4608} + \frac{20161x^4}{82,575,360} - \frac{x^6}{36,700,160} + \frac{3x^8}{587,202,560}\right] +$$

$$\frac{x}{8\sqrt{\pi}} \cdot \exp\left[\frac{-169x^2}{256}\right]\left[1 - \frac{x^2}{4608} + \frac{2629x^4}{9,175,040} - \frac{169x^6}{5,284,823,040} + \frac{28,561x^8}{4,058,744,094,720}\right] +$$

$$\frac{x}{8\sqrt{\pi}} \cdot \exp\left[\frac{-49x^2}{64}\right]\left[1 - \frac{x^2}{4608} + \frac{3049x^4}{9,175,040} - \frac{7x^6}{188,743,680} + \frac{343x^8}{36,238,786,560}\right] +$$

$$\frac{x}{8\sqrt{\pi}} \cdot \exp\left[\frac{-225x^2}{256}\right]\left[1 - \frac{x^2}{4608} + \frac{31501x^4}{82,575,360} - \frac{5x^6}{117,440,512} + \frac{375x^8}{30,064,771,072}\right] +$$

$$\frac{x}{16\sqrt{\pi}} \cdot \exp(-x^2)\left[1 + \frac{127x^2}{4608} + \frac{3929x^4}{9,175,040} + \frac{79x^6}{20,643,840} + \frac{x^8}{61,931,520}\right]$$

When this approximation is utilized with the transition point $x_o = 7.1544$, the relative error bound in the approximation to the error function, over the interval $(0, \infty)$, is $4.82 \times 10^{-16}$.

## Appendix 5:    Proof of Theorem 6.1

Consider the differential equation

$$y_n'(t) + g_n(t)y_n(t) = g_n(t), \qquad y_n(0) = 0, \qquad (159)$$

for the case where $g_n$ is based on the $n$th order approximation $f_n$ to the error function, defined in Theorem 2.1, and is defined according to

$$g_n(t) = \frac{4}{\sqrt{\pi}} e^{-t^2} f_n(t) = \frac{8e^{-t^2}}{\pi} \cdot \sum_{k=0}^{n} c_{n,k} t^{k+1} \left[ p(k, 0) + (-1)^k p(k, t) e^{-t^2} \right], \qquad t \geq 0. \tag{160}$$

To find a solution to the differential equation for such a driving signal, first note that the solution to the differential equation for the case of $g_n(t) = \frac{4}{\sqrt{\pi}} e^{-t^2} \mathrm{erf}(t)$ is

$$y_n(t) = 1 - \exp[-\mathrm{erf}^2(t)]. \tag{161}$$

Second, with $g_n$ defined by Equation 160, the following signal form

$$y_n(t) = 1 - \exp\left[-\left[p_{n,0} + p_{n,1}(t)e^{-t^2} + p_{n,2}(t)e^{-2t^2}\right]\right], \qquad t \geq 0, \tag{162}$$

has potential as a basis for finding the solutions for the unknown polynomial functions $p_{n,1}$ and $p_{n,2}$ and the unknown constant $p_{n,0}$. With such a form, the initial condition of $y_n(0) = 0$ implies

$$p_{n,0} = -[p_{n,1}(0) + p_{n,2}(0)]. \tag{163}$$

It is the case that

$$y_n'(t) = \left[ p_{n,1}^{(1)}(t)e^{-t^2} - 2t p_{n,1}(t)e^{-t^2} + p_{n,2}^{(1)}(t)e^{-2t^2} - 4t p_{n,2}(t)e^{-2t^2} \right] \cdot [1 - y_n(t)]. \tag{164}$$

Substitution of $y_n(t)$ and $y_n'(t)$ into the differential equation yields

$$\left[ p_{n,1}^{(1)}(t)e^{-t^2} - 2t p_{n,1}(t)e^{-t^2} + p_{n,2}^{(1)}(t)e^{-2t^2} - 4t p_{n,2}(t)e^{-2t^2} \right] \exp\left[-\left[p_{n,0} + p_{n,1}(t)e^{-t^2} + p_{n,2}(t)e^{-2t^2}\right]\right] +$$
$$\frac{4}{\sqrt{\pi}} e^{-t^2} f_n(t) \cdot \left[ 1 - \exp\left[-\left[p_{n,0} + p_{n,1}(t)e^{-t^2} + p_{n,2}(t)e^{-2t^2}\right]\right]\right] = \frac{4}{\sqrt{\pi}} e^{-t^2} f_n(t) \tag{165}$$

which implies

$$p_{n,1}^{(1)}(t)e^{-t^2} - 2t p_{n,1}(t)e^{-t^2} + p_{n,2}^{(1)}(t)e^{-2t^2} - 4t p_{n,2}(t)e^{-2t^2} - \frac{4}{\sqrt{\pi}} e^{-t^2} f_n(t) = 0 \tag{166}$$

and

$$p_{n,1}^{(1)}(t)e^{-t^2} - 2t p_{n,1}(t)e^{-t^2} + p_{n,2}^{(1)}(t)e^{-2t^2} - 4t p_{n,2}(t)e^{-2t^2} =$$
$$\frac{8}{\pi} e^{-t^2} \sum_{k=0}^{n} c_{n,k} t^{k+1} \left[ p(k, 0) + (-1)^k p(k, t) e^{-t^2} \right]. \tag{167}$$

Thus:

$$p_{n,1}^{(1)}(t) - 2t p_{n,1}(t) = \frac{8}{\pi} \sum_{k=0}^{n} c_{n,k} p(k, 0) t^{k+1}$$
$$p_{n,2}^{(1)}(t) - 4t p_{n,2}(t) = \frac{8}{\pi} \sum_{k=0}^{n} c_{n,k} (-1)^k p(k, t) t^{k+1} \tag{168}$$

To solve for the polynomials $p_{n,1}$ and $p_{n,2}$, first note (see Equation 26) that

$$p(k, t) = a_{k,0} + a_{k,1} t + a_{k,2} t^2 + \ldots + a_{k,k} t^k, \qquad k \in \{0, 1, \ldots, n\}, \tag{169}$$

for appropriately defined coefficients $a_{k,j}$, $j \in \{0, 1, \ldots, k\}$.

### A5.1     Solving for Coefficients of First Polynomial

Substitution of $p(k, 0)$ from Equation 169, into the differential equation defining $p_{n,1}$, yields

$$p_{n,1}^{(1)}(t) - 2t p_{n,1}(t) = \frac{8}{\pi} \sum_{k=0}^{n} c_{n,k} a_{k,0} t^{k+1}. \tag{170}$$

With $a_{n,0} = 0$ for $n$ odd, the maximum power for $t$ on the right hand side of the differential equation is $t^{n+1}$, $n$ even, and $t^n$ for $n$ odd. Thus, the form required for $p_{n,1}$ is

$$p_{n,1}(t) = \begin{cases} \alpha_0 + \alpha_1 t + \ldots + \alpha_{n-1} t^{n-1} & n \text{ odd} \\ \alpha_0 + \alpha_1 t + \ldots + \alpha_n t^n & n \text{ even} \end{cases} \tag{171}$$

Substitution then yields

$$[\alpha_1 + 2\alpha_2 t + \ldots + (n-1)\alpha_{n-1} t^{n-2}] - 2t[\alpha_0 + \alpha_1 t + \ldots + \alpha_{n-1} t^{n-1}] = \frac{8}{\pi} \sum_{k=0}^{n-1} c_{n,k} a_{k,0} t^{k+1} \quad n \text{ odd}$$

$$[\alpha_1 + 2\alpha_2 t + \ldots + n\alpha_n t^{n-1}] - 2t[\alpha_0 + \alpha_1 t + \ldots + \alpha_n t^n] = \frac{8}{\pi} \sum_{k=0}^{n} c_{n,k} a_{k,0} t^{k+1} \quad n \text{ even} \tag{172}$$

For the case of $n$ even, equating coefficients associated with set powers of $t$, yields:

$$\begin{array}{lll}
t^{n+1}: & -2\alpha_n = \frac{8}{\pi} \cdot c_{n,n} a_{n,0} & \Rightarrow \quad \alpha_n = \frac{-4}{\pi} \cdot c_{n,n} a_{n,0} \\
t^n: & -2\alpha_{n-1} = \frac{8}{\pi} \cdot c_{n,n-1} a_{n-1,0} & \Rightarrow \quad \alpha_{n-1} = \frac{-4}{\pi} \cdot c_{n,n-1} a_{n-1,0} \\
t^{n-1}: & n\alpha_n - 2\alpha_{n-2} = \frac{8}{\pi} \cdot c_{n,n-2} a_{n-2,0} & \Rightarrow \quad \alpha_{n-2} = \frac{n}{2} \cdot \alpha_n - \frac{4}{\pi} \cdot c_{n,n-2} a_{n-2,0} \\
\ldots & & \\
t^2: & 3\alpha_3 - 2\alpha_1 = \frac{8}{\pi} \cdot c_{n,1} a_{1,0} & \Rightarrow \quad \alpha_1 = \frac{3\alpha_3}{2} - \frac{4}{\pi} \cdot c_{n,1} a_{1,0} \\
t^1: & 2\alpha_2 - 2\alpha_0 = \frac{8}{\pi} \cdot c_{n,0} a_{0,0} & \Rightarrow \quad \alpha_0 = \alpha_2 - \frac{4}{\pi} \cdot c_{n,0} a_{0,0}
\end{array} \tag{173}$$

With the odd coefficients $a_{1,0}, a_{3,0}, \ldots, a_{n-1,0}$ being zero, it follows that the corresponding odd coefficients $\alpha_{n-1}, \alpha_{n-3}, \ldots, \alpha_1$ are also zero. For the even coefficients, the algorithm is:

$$\alpha_m = \frac{-4}{\pi} \cdot c_{n,m} a_{m,0}$$

$$\alpha_{m-2i} = \frac{(m-2i+2)\alpha_{m-2i+2}}{2} - \frac{4}{\pi} \cdot c_{n,m-2i} a_{m-2i,0}, \quad i \in \left\{1, \ldots, \frac{m}{2} - 1\right\} \tag{174}$$

$$\alpha_0 = \alpha_2 - \frac{4}{\pi} \cdot c_{n,0} a_{0,0}$$

where $m = n$. For the case of $n$ being odd, the odd coefficients $\alpha_n, \alpha_{n-2}, \ldots, \alpha_1$ are again zero and the algorithm is the same as that specified in Equation 174 with $m = n-1$.

### A5.2     Solving for Coefficients of Second Polynomial

Substitution of $p(k, t)$ from Equation 169, into the differential equation defining $p_{n,2}$, yields:

$$p_{n,2}^{(1)}(t) - 4tp_{n,2}(t) = \frac{8}{\pi} \cdot \sum_{k=0}^{n} c_{n,k}(-1)^k \left[ \sum_{i=0}^{k} a_{k,i} t^{k+i+1} \right]. \tag{175}$$

The coefficients $a_{k,i}$ that are associated with a given power of $t$ are illustrated in Figure 28. It then follows, for a fixed power of $t$, say $t^r$, that the associated coefficients, $a_{k,i}$, are

$$k \in \left\{ \left\lfloor \frac{r}{2} \right\rfloor, \left\lfloor \frac{r}{2} \right\rfloor + 1, \ldots, \min\{r-1, n\} \right\} \tag{176}$$

$$i = r - k - 1.$$

Thus:

$$p_{n,2}^{(1)}(t) - 4tp_{n,2}(t) = \frac{8}{\pi} \cdot \sum_{r=1}^{2n+1} \left[ \sum_{k=\lfloor r/2 \rfloor}^{\min\{r-1,n\}} c_{n,k}(-1)^k a_{k,r-k-1} \right] t^r. \tag{177}$$

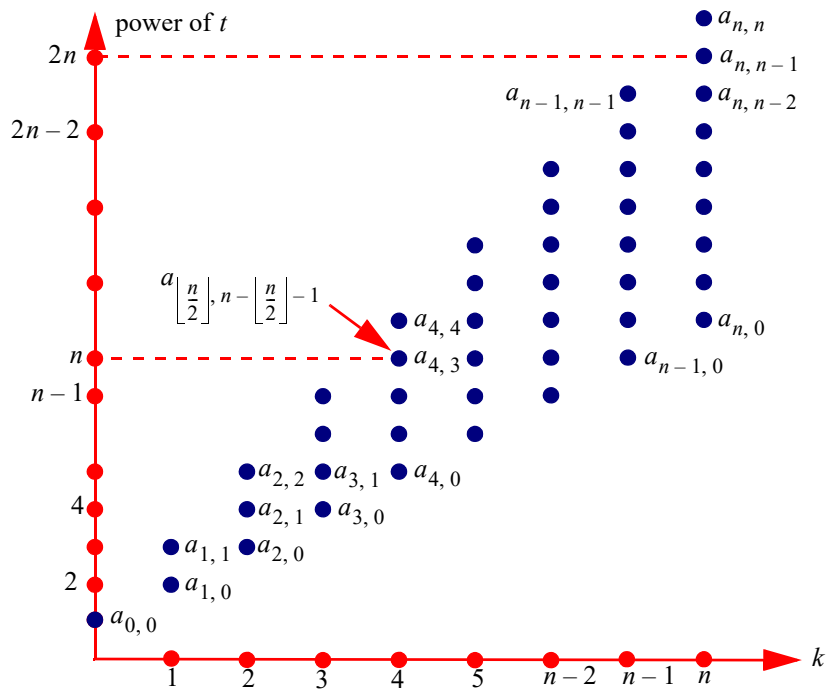

Figure 28.  Illustration of the coefficients that potentially are non-zero for a set power of $t$. The illustration is for the case of $n = 8$.

With

$$p_{n,2}(t) = \beta_0 + \beta_1 t + \ldots + \beta_m t^m \tag{178}$$

the differential equation implies

$$[\beta_1 + 2\beta_2 t + \ldots + m\beta_m t^{m-1}] - 4t[\beta_0 + \beta_1 t + \ldots + \beta_m t^m]$$
$$= \frac{8}{\pi} \cdot \sum_{r=1}^{2n+1} \left[ \sum_{k=\lfloor r/2 \rfloor}^{\min\{r-1,n\}} c_{n,k}(-1)^k a_{k,r-k-1} \right] t^r. \tag{179}$$

The requirement, thus, is for $m = 2n$. Equating coefficients (see Figure 28) yields:

$$t^{2n+1}: \quad -4\beta_{2n} = \frac{8}{\pi} \cdot (-1)^n c_{n,n} a_{n,n} \quad \Rightarrow \quad \beta_{2n} = \frac{-2}{\pi} \cdot (-1)^n c_{n,n} a_{n,n}$$

$$t^{2n}: \quad -4\beta_{2n-1} = \frac{8}{\pi} \cdot [(-1)^n c_{n,n} a_{n,n-1}] \quad \Rightarrow \quad \beta_{2n-1} = \frac{-2}{\pi} \cdot (-1)^n c_{n,n} a_{n,n-1}$$

$$t^{2n-1}: \quad 2n\beta_{2n} - 4\beta_{2n-2} = \frac{8}{\pi} \cdot [(-1)^{n-1} c_{n,n-1} a_{n-1,n-1} + (-1)^n c_{n,n} a_{n,n-2}]$$

$$\Rightarrow \quad \beta_{2n-2} = \frac{2n\beta_{2n}}{4} - \frac{2}{\pi} \cdot [(-1)^{n-1} c_{n,n-1} a_{n-1,n-1} + (-1)^n c_{n,n} a_{n,n-2}] \quad (180)$$

...

$$t^3: \quad 4\beta_4 - 4\beta_2 = \frac{8}{\pi} \cdot [-c_{n,1} a_{1,1} + c_{n,2} a_{2,0}] \quad \Rightarrow \quad \beta_2 = \beta_4 - \frac{2}{\pi} \cdot [-c_{n,1} a_{1,1} + c_{n,2} a_{2,0}]$$

$$t^2: \quad 3\beta_3 - 4\beta_1 = \frac{-8}{\pi} \cdot c_{n,1} a_{1,0} \quad \Rightarrow \quad \beta_1 = \frac{3\beta_3}{4} + \frac{2}{\pi} \cdot c_{n,1} a_{1,0}$$

$$t^1: \quad 2\beta_2 - 4\beta_0 = \frac{8}{\pi} \cdot c_{n,0} a_{0,0} \quad \Rightarrow \quad \beta_0 = \frac{2\beta_2}{4} - \frac{2}{\pi} \cdot c_{n,0} a_{0,0}$$

As the coefficients $a_{n,n-1}, a_{n,n-3}, \ldots, a_{1,0}$ are zero, the algorithm is:

$$\beta_{2n} = \frac{-2}{\pi} \cdot (-1)^n c_{n,n} a_{n,n}$$

$$\beta_{2n-2i} = \frac{[2n-2i+2]\beta_{2n-2i+2}}{4} - \frac{2}{\pi} \sum_{k=n-i}^{\min\{2n-2i,\,n\}} (-1)^k c_{n,k} a_{k,2(n-i)-k}, \quad i \in \{1, \ldots, n-1\} \quad (181)$$

$$\beta_0 = \frac{\beta_2}{2} - \frac{2}{\pi} \cdot c_{n,0} a_{0,0}$$

## Appendix 6:   Proof of Theorem 7.1

Consider the results stated in Theorem 2.1:

$$\text{erf}(x) = \frac{2}{\sqrt{\pi}} \cdot \int_0^x e^{-\lambda^2} d\lambda \approx \frac{2}{\sqrt{\pi}} \cdot \sum_{k=0}^n c_{n,k} x^{k+1} \left[ p(k,0) + (-1)^k p(k,x) e^{-x^2} \right]. \quad (182)$$

Differentiation then yields

$$e^{-x^2} \approx \sum_{k=0}^n c_{n,k}(k+1) x^k \left[ p(k,0) + (-1)^k p(k,x) e^{-x^2} \right] +$$
$$\sum_{k=0}^n c_{n,k}(-1)^k x^{k+1} \left[ p^{(1)}(k,x) e^{-x^2} - 2x p(k,x) e^{-x^2} \right] \quad (183)$$

which leads to the required result:

$$e^{-x^2} \approx \frac{\sum_{k=0}^n c_{n,k}(k+1) x^k p(k,0)}{1 + \sum_{k=0}^n c_{n,k}(-1)^{k+1} x^k [p(k,x)[k+1-2x^2)] + x p^{(1)}(k,x)]} \quad (184)$$

## Appendix 7:   Proof of Theorem 7.2

Consider the exact result

$$\text{erf}(x) = f_0(x) + \varepsilon_0(x) \tag{185}$$

where $f_0$ is specified by Equation 32 and the derivative of the error term, $\varepsilon_0^{(1)}(x)$, is specified by Equation 40. By utilizing a Taylor series approximation for $\exp(-x^2)$, $\varepsilon_0^{(1)}(x)$ can be written as

$$\varepsilon_0^{(1)}(x) = \frac{x^2}{\sqrt{\pi}} \cdot \left[1 - \frac{3x^2}{2} + \frac{5x^4}{6} - \frac{7x^6}{24} + \frac{9x^8}{120} - \ldots + \frac{(-1)^k(2k+1)x^{2k}}{(k+1)!} + \ldots\right] \tag{186}$$

Integration yields

$$\varepsilon_0(x) = \frac{x^3}{\sqrt{\pi}} \cdot \left[\frac{1}{3} - \frac{3x^2}{2 \cdot 5} + \frac{5x^4}{6 \cdot 7} - \frac{7x^6}{24 \cdot 9} + \frac{9x^8}{120 \cdot 11} - \ldots + \frac{(-1)^k(2k+1)x^{2k}}{(2k+3)(k+1)!} + \ldots\right] \tag{187}$$

and the following series for the error function then follows:

$$\begin{aligned}\text{erf}(x) = {}& \frac{x}{\sqrt{\pi}} + \frac{x}{\sqrt{\pi}} \cdot e^{-x^2} + \\ & \frac{x^3}{\sqrt{\pi}} \cdot \left[\frac{1}{3} - \frac{3x^2}{2 \cdot 5} + \frac{5x^4}{6 \cdot 7} - \frac{7x^6}{24 \cdot 9} + \frac{9x^8}{120 \cdot 11} - \ldots + \frac{(-1)^k(2k+1)x^{2k}}{(2k+3)(k+1)!} + \ldots\right]\end{aligned} \tag{188}$$

The series associated with first and second order approximations follow in an analogous manner.

## Appendix 8:   Proof of Theorem 7.3

The filter output is given by the convolution integral:

$$\begin{aligned}y(t) &= \int_0^t \text{erf}\left[\frac{\lambda}{\gamma}\right] \cdot \frac{(t-\lambda)e^{-[t-\lambda]/\tau}}{\tau^2} d\lambda \\ &= \frac{e^{-t/\tau}}{\tau^2} \cdot \left[t\int_0^t \text{erf}\left[\frac{\lambda}{\gamma}\right]e^{\lambda/\tau}d\lambda - \int_0^t \text{erf}\left[\frac{\lambda}{\gamma}\right]\lambda e^{\lambda/\tau}d\lambda\right].\end{aligned} \tag{189}$$

Using the integral results, (e.g. Ng 1969, eqn. 4.2.1 and eqn. 4.2.5)

$$\int_0^t \text{erf}(a\lambda)e^{b\lambda}d\lambda = \frac{1}{b}\text{erf}(at)e^{bt} - \frac{1}{b}\exp\left[\frac{b^2}{4a^2}\right]\text{erf}\left[at - \frac{b}{2a}\right] + \frac{1}{b}\exp\left[\frac{b^2}{4a^2}\right]\text{erf}\left[-\frac{b}{2a}\right] \tag{190}$$

$$\begin{aligned}\int_0^t \text{erf}(a\lambda)\lambda e^{b\lambda}d\lambda = {}& \frac{1}{b}\left[t - \frac{1}{b}\right]\text{erf}(at)e^{bt} - \frac{1}{b}\exp\left[\frac{b^2}{4a^2}\right]\left[\left[\frac{b}{-2a^2} - \frac{1}{b}\right]\text{erf}\left[at - \frac{b}{2a}\right] - \frac{1}{a\sqrt{\pi}}\exp\left[-\left[at - \frac{b}{2a}\right]^2\right]\right] + \\ & \frac{1}{b}\exp\left[\frac{b^2}{4a^2}\right]\left[\left[\frac{b}{-2a^2} - \frac{1}{b}\right]\text{erf}\left[-\frac{b}{2a}\right] - \frac{1}{a\sqrt{\pi}}\exp\left[\frac{-b^2}{4a^2}\right]\right]\end{aligned} \tag{191}$$

with $a = 1/\gamma$, $b = 1/\tau$, $b^2/4a^2 = \gamma^2/4\tau^2$, it then follows that

$$\begin{aligned}y(t) = {}& \frac{te^{-t/\tau}}{\tau^2}\left[\tau\text{erf}\left[\frac{t}{\gamma}\right]e^{t/\tau} - \tau\exp\left[\frac{\gamma^2}{4\tau^2}\right]\text{erf}\left[\frac{t}{\gamma} - \frac{\gamma}{2\tau}\right] + \tau\exp\left[\frac{\gamma^2}{4\tau^2}\right]\text{erf}\left[\frac{-\gamma}{2\tau}\right]\right] - \\ & \frac{e^{-t/\tau}}{\tau^2}\left[\begin{aligned}&\tau(t-\tau)\text{erf}\left[\frac{t}{\gamma}\right]e^{t/\tau} - \tau\exp\left[\frac{\gamma^2}{4\tau^2}\right]\left[\left[\frac{\gamma^2}{2\tau} - \tau\right]\text{erf}\left[\frac{t}{\gamma} - \frac{\gamma}{2\tau}\right] - \frac{\gamma}{\sqrt{\pi}}\exp\left[-\left[\frac{t}{\gamma} - \frac{\gamma}{2\tau}\right]^2\right]\right] + \\ & \tau\exp\left[\frac{\gamma^2}{4\tau^2}\right]\left[\left[\frac{\gamma^2}{2\tau} - \tau\right]\text{erf}\left[\frac{-\gamma}{2\tau}\right] - \frac{\gamma}{\sqrt{\pi}}\exp\left[\frac{-\gamma^2}{4\tau^2}\right]\right]\end{aligned}\right]\end{aligned} \tag{192}$$

Simplifying, and using the fact that the error function is an odd function, yields the required result:

$$y(t) = \text{erf}\left[\frac{t}{\gamma}\right] + \frac{e^{-t/\tau}}{\tau} \cdot \left[\exp\left[\frac{\gamma^2}{4\tau^2}\right]\left[\frac{\gamma^2}{2\tau} - (t+\tau)\right]\left[\text{erf}\left[\frac{\gamma}{2\tau}\right] - \text{erf}\left[\frac{\gamma}{2\tau} - \frac{t}{\gamma}\right]\right] - \frac{\gamma}{\sqrt{\pi}}\exp\left[\frac{\gamma^2}{4\tau^2}\right]\exp\left[-\left[\frac{t}{\gamma} - \frac{\gamma}{2\tau}\right]^2\right] + \frac{\gamma}{\sqrt{\pi}}\right] \tag{193}$$

To prove convergence, consider

$$\lim_{n \to \infty} y_n(t) = \lim_{n \to \infty} \int_0^t f_n\left[\frac{\lambda}{\gamma}\right] h(t-\lambda)d\lambda = \int_0^t \text{erf}\left[\frac{\lambda}{\gamma}\right] h(t-\lambda)d\lambda \tag{194}$$

where $\lim_{n \to \infty} f_n(x) = \text{erf}(x)$ and $h$ is the impulse response of the second order filter. The interchange of limit and integration is valid, consistent with Lemma 2, as the integrand comprises of differentiable bounded functions.